\documentclass[10pt]{article}
\usepackage{amssymb}

\usepackage{amsmath}

\newtheorem{theorem}{Theorem}[section]

\newtheorem{definition}[theorem]{Definition}

\newtheorem{lemma}[theorem]{Lemma}

\newtheorem{proposition}[theorem]{Proposition}

\newenvironment{proof}[1][Proof]{\textbf{#1.} }{\ \rule{0.5em}{0.5em}}

\begin{document}

\title{Stochastic Heat Equation Driven by Fractional Noise and Local Time}
\author{Yaozhong Hu\thanks{%
Y. Hu is supported by the National Science Foundation under DMS0504783 }%
\quad and\quad Nualart\thanks{%
D. Nualart is supported by the National Science Foundation under DMS0604207 }
\\
Department of Mathematics\thinspace ,\ University of Kansas\\
405 Snow Hall\thinspace ,\ Lawrence, Kansas 66045-2142\\
hu@math.ku.edu and nualart@math.ku.edu}
\date{}
\maketitle

\begin{abstract}
The aim of this paper is to study the $d$-dimensional stochastic heat
equation with a multiplicative Gaussian noise which is white in space and it
has the covariance of a fractional Brownian motion with Hurst parameter $%
H\in (0,1)$ in time. Two types of equations are considered. First we
consider the equation in the It\^{o}-Skorohod sense, and later in the
Stratonovich sense. An explicit chaos development for the solution is
obtained. On the other hand, the moments of the solution are expressed in
terms of the exponential moments of some weighted intersection local time of
the Brownian motion.
\end{abstract}

\section{Introduction}

This paper deals with the $d$-dimensional stochastic heat equation%
\begin{equation}
\frac{\partial u}{\partial t}=\frac{1}{2}\Delta u+\ u\diamond \frac{%
\partial ^{2}W^{H}}{\partial t\partial x}  \label{01}
\end{equation}%
driven by a Gaussian noise $W^{H}$ which is a white noise in the
spatial variable and a fractional Brownian motion with Hurst
parameter $H\in (0,1)$ in the time variable (see (\ref{z1}) in the
next section for a precise definition of this noise). The initial
condition $u_{0}$ is a bounded continuous function on
$\mathbb{R}^{d}$, and the solution will be a random field
$\{u_{t,x},t\geq 0,x\in \mathbb{R}^{d}\}$. The symbol $\diamond  $
in Equation (\ref{01}) denotes the Wick product. For $H=\frac{1}{2}$, $\frac{%
\partial ^{2}W^{H}}{\partial t\partial x}$ is a space-time white noise, and
in this case, Equation (\ref{01}) coincides with the stochastic heat
equation considered by Walsh (see \cite{Wa}). We know that in this case the
solution exists only in dimension one ($d=1$).

There has been some recent interest in studying stochastic partial
differential equations driven by a fractional noise. Linear stochastic
evolution equations in a Hilbert space driven by an additive cylindrical fBm
with Hurst parameter $H$ were studied by Duncan et al. in \cite{DMD} in the
case $H\in (\frac{1}{2},1)$ and by Tindel et al. in \cite{TTV} in the
general case, where they provide necessary and sufficient conditions for the
existence and uniqueness of an evolution solution. \ In particular, the heat
equation%
\begin{equation*}
\frac{\partial u}{\partial t}=\frac{1}{2}\ \Delta u+\frac{\partial ^{2}W^{H}%
}{\partial t\partial x}
\end{equation*}%
on $\mathbb{R}^{d}$ has a unique solution if and only if $H>\frac{d}{4}$.
The same result holds when one adds to the above equation a nonlinearity of
the form $b(t,x,u)$, where $b$ satisfies the usual linear growth and
Lipschitz conditions in the variable $u$, uniformly with respect to $(t,x)$
(see Maslowski and Nualart in \cite{MN}). The stochastic heat equation on $%
[0,\infty )\times \mathbb{R}^{d}$ with a multiplicative fractional white
noise of Hurst parameter $H=(H_{0},H_{1},\ldots ,H_{d})$ has been studied by
Hu in \cite{Hu} under the conditions $\frac{1}{2}<H_{i}<1$ for $i=0,\ldots ,d
$ and $\sum_{i=0}^{d}H_{i}<d-\frac{2}{2H_{0}-1}$.

The main purpose of this paper is to find conditions on $H$ and $d$ for the
solution to Equation \ (\ref{01}) to exist as a real-valued stochastic
process, and to relate the moments of the solution to the exponential
moments of weighted intersection local times. This relation is based on
Feynman-Kac's formula applied to a regularization of Equation (\ref{01}). In
order to illustrate this fact, consider the particular case $d=1$ and $H=%
\frac{1}{2}$. It is known that there is no Feynman-Kac's formula for the
solution of the one-dimensional stochastic heat equation driven by a
space-time white noise. Nevertheless, using an approximation of the solution
by regularizing the noise we can establish the following formula for the
moments:%
\begin{equation}
E\left[ u_{t,x}^{k}\right] =E^{B}\left[ \prod%
\limits_{j=1}^{k}u_{0}(x+B_{t}^{j})\exp \left( \sum_{i,j=1,i<j}^{k}\
\int_{0}^{t}\ \delta _{0}(B_{s}^{i}-B_{s}^{j})ds\right) \right] ,  \label{d1}
\end{equation}%
for all $k\geq 2$, where $B_{t}$ is a $k$-dimensional Brownian motion
independent of the spaced-time white noise $W^{\frac{1}{2}}$. \ \ In the
case $H>\frac{1}{2}$ and $d\geq 1$, a similar formula \ holds but $%
\int_{0}^{t}\ \delta _{0}(B_{s}^{i}-B_{s}^{j})ds$ \ has to be replaced by
the weighted intersection local time
\begin{equation}
L_{t}=H(2H-1)\int_{0}^{t}\ \int_{0}^{t}\left| s-r\right| ^{2H-2}\delta
_{0}(B_{s}^{i}-B_{r}^{j})dsdt,  \label{lt}
\end{equation}%
where $\left\{ B^{j},j\geq 1\right\} $ are independent $d$-dimensional
Brownian motions (see Theorem \ref{t1}).

The solution of Equation (\ref{01}) has a formal Wiener chaos expansion $%
u_{t,x}=\sum_{n=0}^{\infty }I_{n}(f_{n}(\cdot ,t,x))$.   Then, for the
existence of a real-valued square integrable solution \ we need%
\begin{equation}
\sum_{n=0}^{\infty }n!\left\| f_{n}(\cdot ,t,x)\right\| _{\mathcal{H}%
_{d}^{\otimes n}}^{2}<\infty ,  \label{z2a}
\end{equation}%
where $\mathcal{H}_{d}$ is the Hilbert space associated with the covariance
of the noise $W^{H}$.   It turns out that, if $H>\frac{1}{2}$,  the
asymptotic behavior of the  norms $\left\| f_{n}(\cdot ,t,x)\right\| _{%
\mathcal{H}_{d}^{\otimes n}}$ is similar to the behavior of the \ $n$th
moment of the random variable $L_{t}$ defined in (\ref{lt}). More precisely,
if $u_{0}$ is a constant $K$, for all $n\geq 1$ we have%
\begin{equation*}
\left( n!\right) ^{2}\left\| f_{n}(\cdot ,t,x)\right\| _{\mathcal{H}%
_{d}^{\otimes n}}^{2}=K^{2}E(L_{t}^{n}).
\end{equation*}
These facts leads to the following results:

\begin{itemize}
\item[i)] If $d=1$ and $H>\frac{1}{2}$, the series (\ref{z2a})\ converges,
and there exists a solution to Equation (\ref{01}) which has moments of all
orders that can be expressed in terms of the exponential moments of the
weighted intersection local times $L_t$. In the case $H=\frac{1}{2}$ we just need
the local time of \ a one-dimensional standard Brownian motion (see (\ref{d1}%
)).

\item[ii)] If $H>\frac{1}{2}$ and $d<4H$, \ the norms $\left\| f_{n}(\cdot
,t,x)\right\| _{\mathcal{H}_{d}^{\otimes n}}$ are finite and $%
E(L_{t}^{n})<\infty $ for all $n$. In the particular case $d=2$, \ the
series (\ref{z2a}) converges if $t$ is small enough, and the solution exists
in a small time interval. Similarly, if $d=2$ the random variable $L_{t}$
satisfies $E(\exp \lambda L_{t})<\infty $ if $\lambda $ and $t$ are small
enough.

\item[iii)] If $d=1$ and $\frac{3}{8}<H<\frac{1}{2}$,   the norms $\left\|
f_{n}(\cdot ,t,x)\right\| _{\mathcal{H}_{d}^{\otimes n}}$ are finite and $%
E(L_{t}^{n})<\infty $ for all $n$.
\end{itemize}

A natural problem is to investigate what happens if we replace the
Wick product by the ordinary product in Equation (\ref{01}), that
is, we  consider the equation

\begin{equation}
\frac{\partial u}{\partial t}=\frac{1}{2}\Delta u+\ u\frac{\partial ^{2}W^{H}%
}{\partial t\partial x}.  \label{02}
\end{equation}%
In terms of the mild formulation, the Wick product leads to the use of It%
\^{o}-Skorohod stochastic integrals, whereas the ordinary product requires
the use of Stratonovich integrals. For this reason, if we use the ordinary
product we must assume $d=1$ and $H>\frac{1}{2}$. In this case we show that
the solution exists and its moments can be computed in terms of exponential
moments of weighted intersection local times and weighted self-intersection
local times in the case $H>\frac{3}{4}$.

The paper is organized as follows. Section 2 contains some preliminaries on
the fractional noise $W^{H}$ and the Skorohod integral with respect to it.
In Section 3 we present the results on the moments of the weighted
intersection local times assuming $H\geq \frac{1}{2}$. \ Section 4 is
devoted to study the Wiener chaos expansion of the solution to Equation \ (%
\ref{01}).   The case $H<\frac{1}{2}$ is more involved because it requires
the use of fractional derivatives. We show here that if $\frac 38 < H< \frac 12$, the norms
$ \| f_{n}(\cdot ,t,x) \| _{\mathcal{H}_{d}^{\otimes n}}$ are finite and they are related to the moments of a fractional derivative of the intersection local time.  We
derive the formulas for the moments of the solution in the case $H\geq \frac{%
1}{2}$ in Section 5. Finally, Section 6 deals with equations defined using
ordinary product and Stratonovich integrals.

\setcounter{equation}{0}
\section{Preliminaries}

Suppose that $W^{H}=\{W^{H}(t,A),t\geq 0,A\in \mathcal{B}(\mathbb{R}%
^{d}),|A|<\infty \}$, where $\mathcal{B}(\mathbb{R}^{d})$ is the Borel $%
\sigma $-algebra of $\mathbb{R}^{d}$, is a zero mean Gaussian family of
random variables with the covariance function%
\begin{equation}
E(W^{H}(t,A)W^{H}(s,B))=\frac{1}{2}(t^{2H}+s^{2H}-|t-s|^{2H})|A\cap B|,
\label{z1}
\end{equation}%
defined in a complete probability space $(\Omega ,\mathcal{F},P)$,
where $H\in (0,1)$, \ and $|A|$ denotes the Lebesgue measure of $A$. Thus, \
for each Borel set $A$ with finite Lebesgue measure, $\{W^{H}(t,A),t\geq 0\}$
is a fractional Brownian motion (fBm) with Hurst parameter $H$ and variance
\ $t^{2H}|A|$, and the fractional Brownian motions corresponding to disjoint
sets are independent.

Then, the multiplicative noise $\frac{\partial ^{2}W^{H}}{\partial t\partial
x}$ appearing in Equation (\ref{01}) is the formal derivative of the random
measure $W^{H}(t,A)$:%
\begin{equation*}
W^{H}(t,A)=\int_{A}\int_{0}^{t}\frac{\partial ^{2}W^{H}}{\partial s\partial x%
}dsdx.
\end{equation*}%
We know that there is an integral representation of the form%
\begin{equation*}
W^{H}(t,A)=\int_{0}^{t}\int_{A}K_{H}(t,s)W(ds,dx),
\end{equation*}%
where $W$ is a space-time white noise, and the square integrable kernel $%
K_{H}$ is given by
\begin{equation*}
K_{H}(t,s)=c_{H}s^{\frac{1}{2}-H}\int_{s}^{t}(u-s)^{H-\frac{3}{2}}u^{H-\frac{%
1}{2}}du,
\end{equation*}%
for some constant $c_{H}$. We will set $K_{H}(t,s)=0$ if $s>t$.

Denote by $\mathcal{E}$ the space of step functions on $\mathbb{R}_{+}$. Let
$\mathcal{H}$ be the closure of $\mathcal{E}$ with respect to the inner
product induced by%
\begin{equation*}
\left\langle \mathbf{1}_{[0,t]},\mathbf{1}_{[0,s]}\right\rangle _{\mathcal{H}%
}=K_{H}(t,s).
\end{equation*}%
The operator $K_{H}^{\ast }:\mathcal{E}\rightarrow L^{2}(\mathbb{R}_{+})$
defined by $K_{H}^{\ast }(\mathbf{1}_{[0,t]})(s)=K_{H}(t,s)$ $\ $provides a
linear isometry between $\mathcal{H}$ and $L^{2}(\mathbb{R}_{+})$.

The mapping $\mathbf{1}_{[0,t]\times A}\rightarrow W^{H}(t,A)$ extends to a
linear isometry between the tensor product $\mathcal{H}\otimes L^{2}(\mathbb{%
R}^{d})$, denoted by $\mathcal{H}_{d}$, and the Gaussian space spanned by $%
W^{H}$. \ We will denote this isometry by $W^{H}$. Then, for each $\varphi
\in \mathcal{H}_{d}$ \ we have%
\begin{equation*}
W^{H}(\varphi )=\int_{0}^{\infty }\int_{\mathbb{R}^{d}}\left( K_{H}^{\ast
}\otimes I\right) \varphi (t,x)W(dt,dx).
\end{equation*}%
We will make use of the notation $W^{H}(\varphi )=\int_{0}^{\infty }\int_{%
\mathbb{R}^{d}}\varphi dW^{H}$.

If $H=\frac{1}{2}$, \ then $\mathcal{H}=L^{2}(\mathbb{R}_{+})$, and the
operator $K_{H}^{\ast }$ is the identity. In this case, we have $\mathcal{H}%
_{d}=L^{2}(\mathbb{R}_{+}\times \mathbb{R}^{d})$. \

Suppose now that $H>\frac{1}{2}$. The operator $K_{H}^{\ast }$ can be
expressed as a fractional integral operator composed with power functions
(see \cite{Nu}). More precisely, for any function $\varphi \in \mathcal{E}$
with support included in the time interval $[0,T]$ we have%
\begin{equation*}
\left( K_{H}^{\ast }\varphi \right) (t)=c_{H}^{\prime }t^{\frac{1}{2}%
-H}I_{T-}^{H-\frac{1}{2}}\left( \varphi (s)s^{H-\frac{1}{2}}\right) (t),
\end{equation*}%
where $I_{T-}^{H-\frac{1}{2}}$ is the right-sided fractional integral
operator defined by%
\begin{equation*}
I_{T-}^{H-\frac{1}{2}}f(t)=\frac{1}{^{\Gamma (H-\frac{1}{2})}}%
\int_{t}^{T}(s-t)^{H-\frac{3}{2}}f(s)ds.
\end{equation*}%
In this case the space $\mathcal{H}$ is not a space of functions (see \cite%
{PT1}) because it contains distributions.   Denote by $|\mathcal{H}|$ the
space of measurable functions  on $[0,T]$ such that
\begin{equation*}
\int_{0}^{\infty }\int_{0}^{\infty }|r-u|^{2H-2}|\varphi _{r}||\varphi
_{u}|drdu<\infty .
\end{equation*}%
Then, $|\mathcal{H}|\subset \mathcal{H}$ and the inner product in the space
\ $\mathcal{H}$ can be expressed in the following form for $\varphi ,\psi
\in |\mathcal{H}|$%
\begin{equation}
\left\langle \varphi ,\psi \right\rangle _{\mathcal{H}}=\int_{0}^{\infty
}\int_{0}^{\infty }\phi (r,u)\varphi _{r}\varphi _{u}drdu,  \label{sc1}
\end{equation}%
where $\phi (s,t)=H(2H-1)|t-s|^{2H-2}$.

Using H\"{o}lder and   Hardy-Littlewood inequalities, \ one can show (see %
\cite{MMV}) that
\begin{equation}
\left\| \varphi \right\| _{\mathcal{H}_{d}}\leq \beta _{H}\left\| \varphi
\right\| _{L^{\frac{1}{H}}(\mathbb{R}_{+};L^{2}(\mathbb{R}^{d}))},
\label{ce3}
\end{equation}%
and this easily implies that$\mathbb{\newline
}$%
\begin{equation}
\left\| \varphi \right\| _{\mathcal{H}_{d}^{\otimes n}}\leq \beta
_{H}^{n}\left\| \varphi \right\| _{L^{\frac{1}{H}}(\mathbb{R}_{+}^{n};L^{2}(%
\mathbb{R}^{nd}))}.  \label{mmv}
\end{equation}

If $H<\frac{1}{2}$, the operator $K_{H}^{\ast }$ can be expressed as a
fractional derivative operator composed with power functions (see \cite{Nu}%
). More precisely, for any function $\varphi \in \mathcal{E}$ with support
included in the time interval $[0,T]$ we have%
\begin{equation*}
\left( K_{H}^{\ast }\varphi \right) (t)=c_{H}^{\prime \prime }t^{\frac{1}{2}%
-H}D_{T-}^{\frac{1}{2}-H}\left( \varphi (s)s^{H-\frac{1}{2}}\right) (t),
\end{equation*}%
where $D_{T-}^{\frac{1}{2}-H}$ is the right-sided fractional derivative
operator defined by%
\begin{equation*}
D_{T-}^{\frac{1}{2}-H}f(t)=\frac{1}{^{\Gamma (H+\frac{1}{2})}}\frac{f(t)}{%
(T-t)^{\frac{1}{2}-H}}-(\frac{1}{2}-H)\int_{t}^{T}\frac{f(s)-f(t)}{(s-t)^{H-%
\frac{3}{2}}}ds.
\end{equation*}%
Moreover, for any $\gamma >\frac{1}{2}-H$ and any $T>0$ we have $C^{\gamma
}([0,T])\subset \mathcal{H}=I_{T-}^{\frac{1}{2}-H}(L^{2}([0,T])$.

If $\varphi $ is a function with support on $[0,T]$, we can express the
operator $K_{H}^{\ast }$ in the following form
\begin{equation}
K_{H}^{\ast }\varphi (t)=K_{H}(T,t)\varphi (t)+\int_{t}^{T}[\varphi
(s)-\varphi (t)]\frac{\partial K_{H}}{\partial s}(s,t)ds.  \label{ke3}
\end{equation}%
We are going to use the following notation for the operator
$K_{H}^{\ast }$:
\begin{equation}
K_{H}^{\ast }\varphi =\int_{[0,T]}\varphi (t)K_{H}^{\ast }(dt,r).
\label{ke1}
\end{equation}%
Notice that if $H>\frac{1}{2}$, the kernel $K_{H}$ vanishes at the diagonal
and we have $\ K_{H}^{\ast }(dt,r)=\frac{\partial K_{H}}{\partial t}(t,r)%
\mathbf{1}_{[r,T]}(t)dt.$

\medskip
Let us now present some preliminaries on the Skorohod integral and the Wick
product. The $n$th Wiener chaos, denoted by $\mathbf{H}_{n}$, is defined as
the closed linear span of the random variables of the form $%
H_{n}(W^{H}(\varphi ))$, where \ $\varphi $ is an element of $\mathcal{H}%
_{d} $ with norm one and $H_{n}$ is the $n$th Hermite polynomial. We denote
by $I_{n}$ the\ linear isometry \ between $\mathcal{H}_{d}^{\otimes n}$
(equipped with the modified norm $\sqrt{n!}\left\| \cdot \right\| _{\mathcal{%
H}_{d}^{\otimes n}}$) and the $n$th Wiener chaos $\mathbf{H}_{n}$, given by $%
I_{n}(\varphi ^{\otimes n})=n!H_{n}(W^{H}(\varphi ))$, for any $\varphi \in
\mathcal{H}_{d}$ with $\left\| \varphi \right\| _{\mathcal{H}_{d}}=1$. Any
square integrable random variable, which is measurable with respect to the $%
\sigma $-field generated by $W^{H}$, has an orthogonal Wiener chaos
expansion of the form%
\begin{equation*}
F=E(F)+\sum_{n=1}^{\infty }I_{n}(f_{n})\text{,}
\end{equation*}%
where $f_{n}$ are symmetric elements of $\mathcal{H}_{d}^{\otimes n}$,
uniquely determined by $F$.

\bigskip Consider a random field $u=\{u_{t,x},t\geq 0,x\in \mathbb{R}^{d}\}$
such that $E\left( u_{t,x}^{2}\right) <\infty $ for all $t,x$. Then, $u$ has
a Wiener chaos expansion of the form%
\begin{equation}
u_{t,x}=E(u_{t,x})+\sum_{n=1}^{\infty }I_{n}(f_{n}(\cdot ,t,x)),  \label{a1}
\end{equation}%
where the series converges in $L^{2}(\Omega )$.

\begin{definition}
We say the random field $u$ satisfying (\ref{a1}) is Skorohod integrable if $%
E(u)\in \mathcal{H}_{d}$, for all $n\geq 1$, \ $f_{n}\in \mathcal{H}%
_{d}^{\otimes (n+1)}$, and the series%
\begin{equation*}
\ W^{H}(E(u))+\sum_{n=1}^{\infty }I_{n+1}(\widetilde{f}_{n})
\end{equation*}%
converges in $L^{2}(\Omega )$, where $\widetilde{f}_{n}$ denotes the
symmetrization of $f_{n}$.We will denote the sum of this series by $\delta
(u)=\int_{0}^{\infty }\int_{\mathbb{R}^{d}}u\delta W^{H}$.
\end{definition}

The Skorohod integral coincides with the adjoint of the derivative operator.
That is, if we define the space $\mathbb{D}^{1,2}$ as  the closure of the
set of smooth and cylindrical random variables of the form%
\begin{equation*}
F=f(W^{H}(h_{1}),\ldots ,W^{H}(h_{n})),
\end{equation*}%
$h_{i}\in \mathcal{H}_{d}$, $f\in C_{p}^{\infty }(\mathbb{R}^{n})$ ($f$ and
all its partial derivatives have polynomial growth) under the norm%
\begin{equation*}
\left\| DF\right\| _{1,2}=\sqrt{E(F^{2})+E(\left\| DF\right\| _{\mathcal{H}%
_{d}}^{2})},
\end{equation*}%
where%
\begin{equation*}
DF=\sum_{j=1}^{n}\frac{\partial f}{\partial x_{j}}(W^{H}(h_{1}),\ldots
,W^{H}(h_{n}))h_{j},
\end{equation*}%
then, the following duality formula holds%
\begin{equation}
E(\delta (u)F)=E\left( \left\langle DF,u\right\rangle _{\mathcal{H}%
_{d}}\right) ,  \label{dua}
\end{equation}%
for any $F\in \mathbb{D}^{1,2}$ and any Skorohod integrable process $u$.

If $F\in \mathbb{D}^{1,2}$ and $h\ $\ is a function which belongs to $%
\mathcal{H}_{d}$, then $Fh$ is Skorohod integrable and, by definition,  the
Wick product equals to the Skorohod integral of $Fh$:%
\begin{equation}
\delta (Fh)=F\diamond  W^{H}(h).  \label{z9}
\end{equation}%
This formula justifies the use of the Wick product in the formulation of
Equation (\ref{01}).

Finally, let us remark that in the case $H=\frac{1}{2}$, if $u_{t,x}$ is an
adapted stochastic process such that $E\left( \int_{0}^{\infty }\int_{%
\mathbb{R}^{d}}u_{t,x}^{2}dxdt\right) <\infty $, then $u$ is Skorohod
integrable and $\delta (u)$ coincides with the It\^{o} stochastic integral:%
\begin{equation*}
\delta (u)=\int_{0}^{\infty }\int_{\mathbb{R}^{d}}u_{t,x}W(dt,dx).
\end{equation*}

\setcounter{equation}{0}
\section{Weighted intersection local times for standard Brownian motions}

In this section we will introduce different kinds of weighted intersection
local times which are relevant in computing the moments of the solutions of
stochastic heat equations with multiplicative fractional noise.

Suppose first that $B^{1}$ and $B^{2}$ are independent \ $d$-dimensional
standard Brownian motions. Consider a nonnegative measurable function $\eta
(s,t)$ on $\mathbb{R}_{+}^{2}$. We are interested in the weighted \
intersection local time formally defined by%
\begin{equation}
I=\int_{0}^{T}\int_{0}^{T}\eta (s,t)\delta _{0}(B_{s}^{1}-B_{t}^{2})dsdt.
\label{x2}
\end{equation}%
We will make use of the following conditions on the weight $\eta $:

\begin{itemize}
\item[C1)] For all $T>0$%
\begin{equation*}
\left\| \eta \right\| _{1,T}:=\max \left( \sup_{0\leq t\leq
T}\int_{0}^{T}\eta (s,t)ds,\sup_{0\leq s\leq T}\int_{0}^{T}\eta
(s,t)dt\right) <\infty .
\end{equation*}

\item[C2)] For all $T>0$ $\ $there exist constants \ $\gamma _{T}>0$ and $%
H\in (0,1)$ such that
\begin{equation*}
\eta (s,t)\leq \gamma _{T}\left| s-t\right| ^{2H-2},
\end{equation*}%
for all $s,t\leq T$.
\end{itemize}

Clearly, C2) is stronger than C1). We will denote by $p_{t}(x)$ the \ $d$%
-dimensional heat kernel $p_{t}(x)=(2\pi t)^{-\frac{d}{2}}e^{-\frac{|x|^{2}}{%
2t}}$. Consider the approximation of \ the intersection local time (\ref{x2}%
) defined by%
\begin{equation}
I_{\varepsilon }=\int_{0}^{T}\int_{0}^{T}\eta (s,t)p_{\varepsilon
}(B_{s}^{1}-B_{t}^{2})dsdt.  \label{a1a}
\end{equation}%
Let us compute the $k$th moment of $I_{\varepsilon }$, where $k\geq 1$ is an
integer. We can write%
\begin{equation}
E\left( I_{\varepsilon }^{k}\right) =\int_{[0,T]^{2k}}\prod_{i=1}^{k}\eta
(s_{i},t_{i})\psi _{\varepsilon }(\mathbf{s},\mathbf{t})d\mathbf{s}d\mathbf{%
t,}  \label{a1b}
\end{equation}%
where $\mathbf{s}=(s_{1},\ldots ,s_{k})$, $\mathbf{t}=(t_{1},\ldots ,t_{k})$
and%
\begin{equation}
\psi _{\varepsilon }(\mathbf{s},\mathbf{t})=E\left( p_{\varepsilon
}(B_{s_{1}}^{1}-B_{t_{1}}^{2})\cdots p_{\varepsilon
}(B_{s_{k}}^{1}-B_{t_{k}}^{2})\right) .  \label{a1c}
\end{equation}%
Using the Fourier transform of the heat kernel we can write%
\begin{eqnarray}
\psi _{\varepsilon }(\mathbf{s},\mathbf{t}) &=&\frac{1}{\left( 2\pi \right)
^{kd}}\left( \int_{\mathbb{R}^{k}}E\left( \exp \left( \sum_{j=1}^{k}\left[
i\left\langle \xi _{j},b_{s_{j}}^{1}-b_{t_{j}}^{2}\right\rangle -\frac{%
\varepsilon }{2}\left| \xi _{j}\right| ^{2}\right] \right) \right) d\mathbf{%
\xi }\right) ^{d}  \notag \\
&=&\frac{1}{\left( 2\pi \right) ^{kd}}\left( \int_{\mathbb{R}^{k}}e^{-\frac{1%
}{2}\sum_{j,l=1}^{k}\xi _{j}\mathrm{Cov}\left(
b_{s_{j}}^{1}-b_{t_{j}}^{2},b_{s_{l}}^{1}-b_{t_{l}}^{2}\right) \xi _{l}-%
\frac{\varepsilon }{2}\left| \mathbf{\xi }\right| ^{2}}d\mathbf{\xi }\right)
^{d},  \label{z3}
\end{eqnarray}%
where $\mathbf{\xi }=(\xi _{1},\ldots ,\xi _{k})$ and $b_{t}^{i}$, $i=1,2$,
are independent one-dimensional Brownian motions. Then $\psi _{\varepsilon }(%
\mathbf{s},\mathbf{t})\leq \psi (\mathbf{s},\mathbf{t})$, where
\begin{equation}
\psi (\mathbf{s},\mathbf{t})=(2\pi )^{-\frac{kd}{2}}\ \left[ \det \left(
s_{j}\wedge s_{l}+t_{j}\wedge t_{l}\right) \right] ^{-\frac{d}{2}}\mathbf{.}
\label{z3a}
\end{equation}%
Set%
\begin{equation}
\alpha _{k}=\int_{[0,T]^{2k}}\prod_{i=1}^{k}\eta (s_{i},t_{i})\psi (\mathbf{s%
},\mathbf{t})d\mathbf{s}d\mathbf{t.}  \label{x7}
\end{equation}%
Then, if $\alpha _{k}<\infty $ for all $k\geq 1$, the family $I_{\varepsilon
}$ converges in $L^{p}$, for all $p\geq 2$, to a limit $I$ and $%
E(I^{k})=\alpha _{k}$. In fact,
\begin{equation*}
\lim_{\varepsilon ,\delta \downarrow 0}E(I_{\varepsilon }I_{\delta })=\alpha
_{2},
\end{equation*}%
so $I_{\varepsilon }$ converges in $L^{2}$, and the convergence in $L^{p}$
follows from the boundedness in $L^{q}$ for $q>p$. Then the following result
holds.

\begin{proposition}
\label{p2} Suppose that C1) holds and $d=1$. Then, for all $\lambda >0$ the
random variable defined in (\ref{a1a}) satisfies%
\begin{equation}
\sup_{\varepsilon >0}E\left( \exp \left( \lambda I_{\varepsilon }\right)
\right) \leq 1+\Phi \left( \frac{\sqrt{T}}{2}\left\| \eta \right\|
_{1,T}\lambda \right) ,  \label{a2}
\end{equation}%
where $\Phi (x)=\sum_{k=1}^{\infty }\frac{x^{k}}{\Gamma (\frac{k+1}{2})}$.
Also, $I_{\varepsilon }$ converges in $L^{p}$ for all $p\geq 2$, and the
limit, denoted by $I$, satisfies the estimate (\ref{a2}).
\end{proposition}

\begin{proof}
The term $\psi (\mathbf{s},\mathbf{t})$ defined in (\ref{z3a}) can be
estimated using  Cauchy-Schwarz inequality:
\begin{eqnarray}
\psi (\mathbf{s},\mathbf{t}) &\leq &(2\pi )^{-k}\left[ \det \left(
s_{j}\wedge s_{l}\right) \right] ^{\frac{1}{4}}\left[ \det \left(
t_{j}\wedge t_{l}\right) \right] ^{\frac{1}{4}}  \notag \\
&=&2^{-k}\pi ^{-\frac{k}{2}}\left[ \beta (\mathbf{s})\beta (\mathbf{t})%
\right] ^{-\frac{1}{4}},  \label{x6a}
\end{eqnarray}%
where for any element $(s_{1},\ldots ,s_{k})\in (0,\infty )^{k}$ with $%
s_{i}\neq s_{j}$ if $i\neq j$, we denote by $\sigma $ the permutation of its
coordinates such that $s_{\sigma (1)}<\cdots <s_{\sigma (n)}$ and $\beta (%
\mathbf{s})=s_{\sigma (1)}(s_{\sigma (2)}-s_{\sigma (1)})\cdots (s_{\sigma
(k)}-s_{\sigma (k-1)})$. Therefore, from (\ref{x6a}) and (\ref{x7}) we obtain%
\begin{equation}
\alpha _{k}\leq 2^{-k}\pi ^{-\frac{k}{2}}\int_{[0,T]^{2k}}\prod_{i=1}^{k}%
\eta (s_{i},t_{i})\left[ \beta (\mathbf{s})\beta (\mathbf{t})\right] ^{-%
\frac{1}{4}}d\mathbf{s}d\mathbf{t.}  \label{x1}
\end{equation}%
Applying again Cauchy-Schwarz inequality yields%
\begin{eqnarray}
\alpha _{k} &\leq &2^{-k}\pi ^{-\frac{k}{2}}\left\{ \int_{[0,T]^{2k}}\
\prod_{i=1}^{k}\eta (s_{i},t_{i})\left[ \beta (\mathbf{s})\right] ^{-\frac{1%
}{2}}d\mathbf{s}d\mathbf{t}\right\} ^{\frac{1}{2}}  \notag \\
&&\times \left\{ \int_{\lbrack 0,T]^{2k}}\ \prod_{i=1}^{k}\eta (s_{i},t_{i})
\left[ \beta (\mathbf{t})\right] ^{-\frac{1}{2}}d\mathbf{s}d\mathbf{t}%
\right\} ^{\frac{1}{2}}  \notag \\
&\leq &\left( 2^{-1}\pi ^{-\frac{1}{2}}\left\| \eta \right\| _{1,T}\right)
^{k}k!\int_{T_{k}}\left[ \beta (\mathbf{s})\right] ^{-\frac{1}{2}}d\mathbf{s}
\notag \\
&=&\frac{k!2^{-k}T^{\frac{k}{2}}\left\| \eta \right\| _{1,T}^{k}}{\Gamma (%
\frac{k+1}{2})},  \label{z6}
\end{eqnarray}%
where $T_{k}=\{\mathbf{s}=(s_{1},\ldots ,s_{k}):0<s_{1}<\cdots <s_{k}<T\}$,
which implies the estimate (\ref{a2}).
\end{proof}

This result can be extended to the case of a $d$-dimensional   Brownian
motion under the stronger condition C2):

\begin{proposition}
\label{prop5} Suppose that C2) holds and $d<4H$. Then, $\lim_{\varepsilon
\downarrow 0}I_{\varepsilon }=I$, exists in $L^{p}$, for all $p\geq 2$. \
Moreover, if $d=2$ and $\lambda <\lambda _{0}(T)$, where
\begin{equation}
\lambda _{0}(T)=\frac{H(2H-1)4\pi }{\gamma _{T}\beta _{H}}T^{\frac{d}{2}%
-2H}\Gamma \left( 1-\frac{d}{4H}\right) ^{-2H},  \label{v0}
\end{equation}%
and $\beta _{H}$ is the constant appearing in the inequality (\ref{ce3}),
then
\begin{equation}
\sup_{\varepsilon >0}E\left( \exp \left( \lambda I_{\varepsilon }\right)
\right) <\infty ,  \label{v1}
\end{equation}%
and $I$ satisfies $E\left( \exp \left( \lambda I\right) \right) <\infty $.
\end{proposition}

\begin{proof}
As in the proof of Proposition \ref{p2}, using condition C2) and inequality (%
\ref{mmv}) we obtain the estimates
\begin{eqnarray*}
\alpha _{k} &\leq &\gamma _{T}^{k}2^{-dk}\pi ^{-\frac{kd}{2}%
}\int_{[0,T]^{2k}}\prod_{i=1}^{k}|t_{i}-s_{i}|^{2H-2}\left[ \beta (\mathbf{s}%
)\beta (\mathbf{t})\right] ^{-\frac{d}{4}}d\mathbf{s}d\mathbf{t} \\
&\leq &\gamma _{T}^{k}2^{-dk}\pi ^{-\frac{kd}{2}}\alpha _{H}^{k}\left(
\int_{[0,T]^{k}}\left[ \beta (\mathbf{s})\right] ^{-\frac{d}{4H}}d\mathbf{s}%
\right) ^{2H} \\
&=&\left( \gamma _{T}\alpha _{H}2^{-2}\pi ^{-1}\right) ^{k}\left( k!\right)
^{2H}\frac{\Gamma (1-\frac{d}{4H})^{k2H}T^{k(1-\frac{d}{4H})2H}}{\Gamma (k(1-%
\frac{d}{4H})+1)^{2H}} \\
&\thickapprox &\left( \gamma _{T}\alpha _{H}2^{-2}\pi ^{-1}\Gamma (1-\frac{d%
}{4H})^{2H}T^{2H-\frac{d}{2}}\right) ^{k}(k!)^{\frac{d}{2}},
\end{eqnarray*}%
where $\alpha _{H}=\frac{\beta _{H}}{H(2H-1)}$. This allows us to conclude
the proof.
\end{proof}

If $d=2$ and $\eta(s,t)=1$ it is   known that the intersection local time
$\int_0^T \int_0^t \delta_0(B^1_s-B^2_t)dsdt$ exists and it has finite  exponential moments up to  a critical exponent  $\lambda_0$ (see Le Gall \cite{Le} and Bass and Chen \cite{BC}).

\medskip
Consider now a one-dimensional standard Brownian motion $B$, and the
weighted self-intersection local time%
\begin{equation*}
I=\int_{0}^{T}\int_{0}^{T}\eta (s,t)\delta _{0}(B_{s}-B_{t})dsdt.
\end{equation*}%
As before, set
\begin{equation*}
I_{\varepsilon }=\int_{0}^{T}\int_{0}^{T}\eta (s,t)p_{\varepsilon
}(B_{s}-B_{t})dsdt.
\end{equation*}

\begin{proposition}
Suppose that C2) holds. {If }$H>\frac{1}{2}$, then we have
\begin{equation}
\sup_{\varepsilon >0}E  \left( \exp \left( {\lambda }\left[ I_{\varepsilon
}-E  \left( I_{\varepsilon }\right) \right] \right) \right) <\infty ,
\label{a4}
\end{equation}%
for all $\lambda >0$. Moreover, the normalized local time $I-E  \left(
I\right) $ exists as a \ limit in $L^{p}$ of $I_{\varepsilon }-E  \left(
I_{\varepsilon }\right) $, for all $p\geq 2$, and it has exponential moments
of all orders.

If $H>\frac{3}{4}$, then we have for all $\lambda >0$
\begin{equation}
\sup_{\varepsilon >0}E  \left( \exp \left( {\lambda }I_{\varepsilon }\right)
\right) <\infty ,  \label{a5}
\end{equation}%
for all $\lambda >0$, and the local time $I$ exists as a \ limit in $L^{p}$
of $I_{\varepsilon }$, for all $p\geq 2$, and it is exponentially integrable.
\end{proposition}

\begin{proof}
We will follow the ideas of Le Gall in \cite{Le}. Suppose first that $H>%
\frac{1}{2}$ and let us show (\ref{a4}).  To simplify the proof we assume $%
T=1$. It suffices to show these results for
\begin{equation*}
J_{\varepsilon }:=\int_{0}^{1}\int_{0}^{t}\eta (s,t)p_{\varepsilon
}(B_{s}-B_{t})dsdt.
\end{equation*}%
Denote, for $n\geq 1$, and $1\leq k\leq 2^{n-1}$
\begin{equation*}
A_{n,k}=\left[ \frac{2k-2}{2^{n}},\frac{2k-1}{2^{n}}\right] \times \left[
\frac{2k-1}{2^{n}},\frac{2k}{2^{n}}\right] .
\end{equation*}%
Set
\begin{equation*}
\alpha _{n,k}^{\varepsilon }=\int_{A_{n,k}}\eta (s,t)p_{\varepsilon
}(B_{s}-B_{t})dsdt
\end{equation*}%
and
\begin{equation*}
\bar{\alpha}_{n,k}^{\varepsilon }=\alpha _{n,k}^{\varepsilon }-E\left(
\alpha _{n,k}^{\varepsilon }\right) .
\end{equation*}%
Notice that the random variables $\alpha _{n,k}^{\varepsilon }$, $1\leq
k\leq 2^{n-1}$, are independent. We have
\begin{equation*}
J_{\varepsilon }=\sum_{n=1}^{\infty }\sum_{k=1}^{2^{n-1}}\alpha
_{n,k}^{\varepsilon },
\end{equation*}%
and%
\begin{equation*}
J_{\varepsilon }-E\left( J_{\varepsilon }\right) =\sum_{n=1}^{\infty
}\sum_{k=1}^{2^{n-1}}\bar{\alpha}_{n,k}^{\varepsilon }.
\end{equation*}%
We can write%
\begin{eqnarray*}
\alpha _{n,k}^{\varepsilon } &=&2^{-2n}\int_{0}^{1}\int_{0}^{1}\eta \left(
\frac{2k-1}{2^{n}}-\frac{s}{2^{n}},\frac{2k-1}{2^{n}}+\frac{t}{2^{n}}\right)  \\
 && \times p_{\varepsilon }(B_{\frac{2k-1}{2^{n}}-\frac{s}{2^{n}}}-B_{\frac{2k-1}{2^{n}}%
+\frac{t}{2^{n}}})dsdt \\
&\leq &\gamma
_{1}2^{-2n-(2H-2)n}\int_{0}^{1}\int_{0}^{1}|t+s|^{2H-2}p_{\varepsilon }(B_{%
\frac{2k-1}{2^{n}}-\frac{s}{2^{n}}}-B_{\frac{2k-1}{2^{n}}+\frac{t}{2^{n}}%
})dsdt,
\end{eqnarray*}%
which has the same distribution as%
\begin{equation*}
\beta _{n,k}^{\varepsilon }=\gamma _{1}2^{-\frac{3}{2}n-(2H-2)n}\int_{0}^{T}%
\int_{0}^{T}|t+s|^{2H-2}p_{\varepsilon 2^{n}}(B_{s}^{1}-B_{t}^{2})dsdt,
\end{equation*}%
where $B^{1}$ and $B^{2}$ are independent one-dimensional Brownian motions.
Hence, using the estimate (\ref{z6}), we obtain%
\begin{eqnarray*}
E\left( \exp \left( {\lambda }\left( \bar{\alpha}_{n,k}^{\varepsilon
}\right) \right) \right)  &=&1+\sum_{j=2}^{\infty }\frac{\lambda ^{j}}{j!}%
E\left( \left( \bar{\alpha}_{n,k}^{\varepsilon }\right) ^{j}\right)  \\
&\leq &1+\sum_{j=2}^{\infty }\frac{\left( 2\lambda \right) ^{j}}{j!}E\left(
\left( \beta _{n,k}^{\varepsilon }\right) ^{j}\right)  \\
&\leq &1+\sum_{j=2}^{\infty }\frac{\left( C_{T}2^{-\frac{3}{2}%
n-(2H-2)n}\lambda \right) ^{j}}{\Gamma (\frac{j+2}{2})},
\end{eqnarray*}%
for some constant $C_{T}$. Hence,%
\begin{equation}
E  \left( \exp \left( {\lambda }\left( \bar{\alpha}_{n,k}^{\varepsilon
}\right) \right) \right) \leq 1+c_{\lambda }2^{-3n-2(2H-2)n}\lambda ^{2},
\label{a6}
\end{equation}%
for some function $c_{\lambda }$.

Fix $a>0$ \ such that $a<2(2H-1\dot{)}$. For any $N\geq 2$ define
\begin{equation*}
b_{N}=\prod_{j=2}^{N}(1-2^{-a(j-1)}),
\end{equation*}%
and notice that $\lim_{N\rightarrow \infty }b_{N}=b_{\infty }>0$. Then, by \
H\"{o}lder's inequality, for all $N\geq 2$ we have
\begin{eqnarray*}
&&E \left[ \exp \left( \lambda b_{N}\sum_{n=1}^{N}\sum_{k=1}^{2^{n-1}}\bar{%
\alpha}_{n,k}^{\varepsilon }\right) \right]  \\
&\leq &\left\{ E \left[ \exp \left( \frac{\lambda b_{N}}{1-2^{-a(N-1)}}%
\sum_{n=1}^{N-1}\sum_{k=1}^{2^{n-1}}\bar{\alpha}_{n,k}^{\varepsilon }\right) %
\right] \right\} ^{1-2^{-a(N-1)}} \\
&&\times \left\{ E \left[ \exp \left( \lambda
b_{N}2^{a(N-1)}\sum_{k=1}^{2^{N-1}}\bar{\alpha}_{N,k}^{\varepsilon }\right) %
\right] \right\} ^{2^{-a(N-1)}} \\
&\leq &\left\{ E \left[ \exp \left( \lambda
b_{N-1}\sum_{n=1}^{N-1}\sum_{k=1}^{2^{n-1}}\bar{\alpha}_{n,k}^{\varepsilon
}\right) \right] \right\}     \\
&&\times \left\{ E \left[ \exp \left( \lambda
b_{N}2^{a(N-1)}\bar{\alpha}_{N,k}^{\varepsilon }\right) \right] \right\}
^{2^{(1-a)(N-1)}}.
\end{eqnarray*}%
Using (\ref{a6}), the second factor in the above expression can be dominated
by
\begin{eqnarray*}
&&\left\{ E \left[ \exp \left( \lambda b_{N}2^{a(N-1)}\bar{\alpha}%
_{N,k}\right) \right] \right\} ^{2^{(1-a)(N-1)}} \\
&\leq &\left( 1+c_{\lambda }\lambda
^{2}b_{2}^{2}2^{2a(N-1)}2^{-3N-2(2H-2)N}\right) ^{2^{(1-a)(N-1)}} \\
&\leq &\exp \left( \kappa c_{\lambda }\lambda ^{2}2^{(a-2-2(2H-2))N}\right) ,
\end{eqnarray*}%
where $\kappa =b_{2}^{2}2^{-a-1}$. Thus by induction we have
\begin{eqnarray*}
E \left[ \exp \left( \lambda b_{N}\sum_{n=1}^{N}\sum_{k=1}^{2^{n-1}}\bar{%
\alpha}_{n,k}\right) \right]  &\leq &\exp \left\{ \sum_{n=2}^{N}\kappa
c_{\lambda }\lambda ^{2}2^{(a-2-2(2H-2))n}\right\}   \\
&&\times E  \left( \exp \bar{\alpha%
}_{1,1}\right)   \\
&\leq &\exp (\kappa c_{\lambda }\lambda ^{2}(1-2^{a+2-4H})^{-1})\\
&& \times E\left( \exp
(\bar{\alpha}_{1,1})\right) <\infty ,
\end{eqnarray*}%
because $\ a<2(2H-1)$. By Fatou lemma we see that
\begin{equation*}
\sup_{\varepsilon >0}E\ \left( \exp \left( \lambda b_{\infty }\left(
J_{\varepsilon }-E\left( J_{\varepsilon }\right) \right) \right) \right)
<\infty ,
\end{equation*}%
and (\ref{a4}) follows.

On the other hand,  one can easily  show that
\begin{eqnarray*}
&& \lim_{\varepsilon ,\delta \downarrow 0}E(\left( J_{\varepsilon }-E\left(
    J_{\varepsilon }\right) \right) \left( J_{\delta }-E\left( J_{\delta
}\right) \right) )
=\frac{1}{2\pi }\int_{s<t<1,s^{\prime }<t^{\prime }<1}\ \eta (s,t)\eta
(s^{\prime },t^{\prime }) \\
&&\times \Bigg[  \left( \det \left[
\begin{array}{cc}
  t-s & \ \left| [s,t]\cap \lbrack s^{\prime },t^{\prime }]\right|  \\
  \left| [s,t]\cap \lbrack s^{\prime },t^{\prime }]\right|  & t^{\prime
}-s^{\prime }%
\end{array}%
\right] \right) ^{-\frac{1}{2}}  \\
&& -\left( (t-s)(t^{\prime }-s^{\prime })\right)
^{-\frac{1}{2}}\Bigg] dsdtds^{\prime }dt^{\prime }
 <\infty ,
\end{eqnarray*}%
which implies the convergence of $I_{\varepsilon }$ in $L^{2}$. The
convergence in $L^{p}$ for $p\geq 2$ and the estimate (\ref{a4}) follow
immediately.

The proof of the inequality (\ref{a5}) is similar. The estimate (\ref{a6})
is replaced by%
\begin{equation}
E \left( \exp \left( {\lambda }\left( \alpha _{n,k}^{\varepsilon }\right)
\right) \right) \leq 1+d_{\lambda }2^{-3n-2(2H-2)n}\lambda ,  \label{a7}
\end{equation}%
for a suitable function $d_{\lambda }$, and we obtain%
\begin{eqnarray*}
&&E \left[ \exp \left( \lambda
b_{N}\sum_{n=1}^{N}\sum_{k=1}^{2^{n-1}}\alpha _{n,k}\right) \right]  \\
&\leq &\exp \left\{ \sum_{n=2}^{N}\sqrt{\kappa }d_{\lambda }\lambda
2^{\left( -\frac{5}{2}-2H\right) n}\right\} E\left( \exp \left( \alpha
_{1,1}\right) \right)  \\
&\leq &\exp (\sqrt{\kappa }d_{\lambda }\lambda ^{2}(1-2^{\left( -\frac{5}{2}%
-2H\right) n})^{-1})E\left( \exp \left( \alpha _{1,1}\right) \right) <\infty
,
\end{eqnarray*}
because $H>\frac{3}{4}$. By Fatou lemma we see that
\begin{equation*}
\sup_{\varepsilon >0}E \left( \exp \left( \lambda b_{\infty }\left(
J_{\varepsilon }-E\left( J_{\varepsilon }\right) \right) \right) \right)
<\infty ,
\end{equation*}%
which implies (\ref{a5}). The convergence in $L^{p}$ of $I_{\varepsilon }$
is proved as usual.
\end{proof}

Notice that condition $H>\frac{3}{4}$ cannot be improved because%
\begin{equation*}
\ E\left( \int_{0}^{T}\int_{0}^{T}|t-s|^{-\frac{1}{2}}\delta
_{0}(B_{s}-B_{t})dsdt\right) =\frac{1}{\sqrt{2\pi }}\int_{0}^{T}%
\int_{0}^{T}|t-s|^{-1}dsdt=\infty .
\end{equation*}

\setcounter{equation}{0}
\section{Stochastic heat equation in the It\^{o}-Skorohod sense}

In this section we study the stochastic partial differential equation (\ref%
{01}) on $\mathbb{R}^{d}$, where $W^{H}$ is a zero mean Gaussian family of
random variables with the covariance function (\ref{z1}), defined on a
complete probability space $(\Omega ,\mathcal{F},P)$, and the initial
condition $u_{0}$ belongs to $C_{b}(\mathbb{R}^{d})$. First we give the
definition of a solution using the Skorhohod integral, which corresponds
formally to the Wick product appearing in Equation (\ref{01}).

For any $t\geq 0$, we denote by $\mathcal{F}_{t}$ the $\sigma $-field
generated by the random variables $\{W(s,A),0\leq s\leq t,A\in \mathcal{B}(%
\mathbb{R}^{d}),|A|<\infty \}$ and the $P$-null sets. A random field $%
u=\{u_{t,x},t\geq 0,x\in \mathbb{R}\}$ is adapted if for any $(t,x)$, $%
u_{t,x}$ is $\mathcal{F}_{t}$-measurable.

For any bounded Borel function $\varphi $ on $\mathbb{R}$ we write $%
p_{t}\varphi (x)=\int_{\mathbb{R}^{d}}p_{t}(x-y)\varphi (y)dy$.

\begin{definition}
\label{def1} An adapted random field $u=\{u_{t,x},t\geq 0,x\in \mathbb{R}%
^{d}\}\,\ $such that $E(u_{t,x}^{2})<\infty \,$ for all $(t,x)$\ is a
solution to Equation (\ref{01}) if \ for any $(t,x)\in \lbrack 0,\infty
)\times \mathbb{R}^{d}$, the process $\{p_{t-s}(x-y)u_{s,y}\mathbf{1}%
_{[0,t]}(s),s\geq 0,y\in \mathbb{R}^{d}\mathbb{\}}$ is Skorohod integrable,
and the \ following equation holds%
\begin{equation}
u_{t,x}=p_{t}u_{0}(x)+\int_{0}^{t}\int_{\mathbb{R}^{d}}p_{t-s}(x-y)u_{s,y}%
\delta W_{s,y}^{H}.  \label{e2}
\end{equation}
\end{definition}

The fact that Equation (\ref{01}) contains a multiplicative Gaussian noise
allows us to find recursively an explicit expression for the Wiener chaos
expansion of the solution. This approach has extensively used in the
literature. For instance, we refer to the papers by Hu \cite{Hu}, Buckdahn
and Nualart \cite{BuNu}, Nualart and Zakai \cite{NuZa}, Nualart and
Rozovskii \cite{NuRo}, \ and Tudor \cite{Tu}, among others.

Suppose that $u=\{u_{t,x},t\geq 0,x\in \mathbb{R}^{d}\mathbb{\}}$ is a
solution to Equation (\ref{01}). Then, for any fixed $(t,x)$, the random
variable $u_{t,x}$ admits the following Wiener chaos expansion
\begin{equation}
u_{t,x}=\sum_{n=0}^{\infty }I_{n}(f_{n}(\cdot ,t,x)),  \label{e3}
\end{equation}%
where for each $(t,x)$, $f_{n}(\cdot ,t,x)$ is a symmetric element in $%
\mathcal{H}_{d}^{\otimes n}$. To find the explicit form of $f_{n}$ we
substitute (\ref{e3}) in the Skorohod integral appearing in (\ref{e2}) we
obtain
\begin{eqnarray*}
\int_{0}^{t}\int_{\mathbb{R}^{d}}p_{t-s}(x-y)u_{s,y}\delta W_{s,y}^{H}
&=&\sum_{n=0}^{\infty }\int_{0}^{t}\int_{\mathbb{R}%
^{r}}I_{n}(p_{t-s}(x-y)f_{n}(\cdot ,s,y))\ \delta W_{s,y}^{H}\  \\
&=&\sum_{n=0}^{\infty }I_{n+1}(\widetilde{p_{t-s}(x-y)f_{n}(\cdot ,s,y)})\ .
\end{eqnarray*}%
Here, $(\widetilde{p_{t-s}(x-y)f_{n}(\cdot ,s,y)}$ denotes the
symmetrization of the function
\begin{equation*}
p_{t-s}(x-y)f_{n}(s_{1},x_{1};\ldots ;s_{n},x_{n};s,y)
\end{equation*}%
in the variables $(s_{1},x_{1}),\ldots ,(s_{n},x_{n}),(s,y)$, that is,
\begin{eqnarray*}
&&\widetilde{p_{t-s}(x-y)f_{n}(\cdot ,s,y)}=\frac{1}{n+1}%
[p_{t-s}(x-y)f_{n}(s_{1},x_{1},\ldots ,s_{n},x_{n},s,y) \\
&&+\sum_{j=1}^{n}p_{t-s_{j}}(x-y_{j})   \\
&&\times f_{n}(s_{1},x_{1},\ldots
,s_{j-1},x_{j-1},s,y,s_{j+1},x_{j+1},\ldots ,s_{n},y_{n},s_{j},y_{j})].
\end{eqnarray*}%
Thus, Equation (\ref{e2}) is equivalent to say that $%
f_{0}(t,x)=p_{t}u_{0}(x) $, and%
\begin{equation}
f_{n+1}(\cdot ,t,x)=\widetilde{p_{t-s}(x-y)f_{n}(\cdot ,s,y)}  \label{e4a}
\end{equation}%
for all $n\geq 0$. Notice that, the adaptability property of the random
field $u$ implies that $f_{n}(s_{1},x_{1},\ldots ,s_{n},x_{n},t,x)=0\,\ $\
if $s_{j}>t$ for some $j$.

This leads to the following formula for the kernels $f_{n}$, for $n\geq 1$%
\begin{eqnarray}
&&f_{n}(s_{1},x_{1},\ldots ,s_{n},x_{n},t,x)=\frac{1}{n!}\   \notag \\
&&\times p_{t-s_{\sigma (n)}}(x-x_{\sigma (n)})\cdots p_{s_{\sigma
(2)}-s_{\sigma (1)}}(x_{\sigma (2)}-x_{\sigma (1)})p_{s_{\sigma
(1)}}u_{0}(x_{\sigma (1)}),  \label{e7}
\end{eqnarray}%
where $\sigma $ denotes the permutation of $\{1,2,\ldots ,n\}$ such that $%
0<s_{\sigma (1)}<\cdots <s_{\sigma (n)}<t$. This implies that there is a
unique solution to Equation (\ref{e2}), and \ the kernels of its chaos
expansion are given by (\ref{e7}). In order to show the existence of a
solution, it suffices to check that the kernels defined in (\ref{e7})
determine an adapted random field satisfying the conditions of Definition %
\ref{def1}. This is equivalent to show that for all $(t,x)$ we have
\begin{equation}
\sum_{n=1}^{\infty }n!\left\| f_{n}(\cdot ,t,x)\right\| _{\mathcal{H}%
_{d}^{\otimes n}}^{2}<\infty .  \label{e5a}
\end{equation}%
It is easy to show that (\ref{e5a}) holds if $H=\frac{1}{2}$ and $d=1$. In
fact, we have, assuming $|u_{0}|\leq K$, and with the notation $\mathbf{x}%
=(x_{1},\ldots ,x_{n})$, and $\mathbf{s}=(s_{1},\ldots ,s_{n})$:
\begin{eqnarray*}
&&\left\| f_{n}(\cdot ,t,x)\right\| _{\mathcal{H}_{1}^{\otimes n}}^{2} \\
&=&\frac{1}{\left( n!\right) ^{2}}\int_{[0,t]^{n}}\int_{\mathbb{R}%
^{n}}p_{t-s_{\sigma (n)}}(x-x_{\sigma (n)})^{2}\cdots p_{s_{\sigma
(2)}-s_{\sigma (1)}}(x_{\sigma (2)}-x_{\sigma (1)})^{2} \\
&&\times p_{s_{\sigma (1)}}u_{0}(x_{\sigma (1)})^{2}\ d\mathbf{x}d\mathbf{s}
\\
&\leq &K^{2}\frac{\left( 4\pi \right) ^{-\frac{n}{2}}}{\left( n!\right) ^{2}}%
\int_{[0,t]^{n}}\prod\limits_{j=1}^{n}(s_{\sigma (j+1)}-s_{\sigma (j)})^{-%
\frac{1}{2}}d\mathbf{s} \\
&=&\ \frac{K^{2}\left( 4\pi \right) ^{-\frac{n}{2}}}{n!}\int_{T_{n}}\prod%
\limits_{j=1}^{n}(s_{j+1}-s_{j})^{-\frac{1}{2}}d\mathbf{s,}
\end{eqnarray*}%
where $T_{n}=\{(s_{1},\ldots ,s_{n})\in \lbrack 0,t]^{n}:0<s_{1}<\cdots
<s_{n}<t\}$ and by convention $s_{n+1}=t$. Hence,
\begin{equation*}
\left\| f_{n}(\cdot ,t,x)\right\| _{\mathcal{H}_{1}^{\otimes n}}^{2}\leq
\frac{K^{2}2^{-n}t^{\frac{n}{2}}}{n!\Gamma (\frac{n+1}{2})},
\end{equation*}%
which implies (\ref{e5a}). On the other hand, if $H=\frac{1}{2}$ and $d\geq
2 $, these norms are infinite.

Notice that if $u_{0}=1$, then \ $(n!)^{2}\left\| f_{n}(\cdot ,t,x)\right\|
_{\mathcal{H}_{1}^{\otimes n}}^{2}$ coincides with the moment of order $n$
of the local time at zero    of the one-dimensional Brownian
motion with variance $2t$, that is,%
\begin{equation*}
(n!)^{2}\left\| f_{n}(\cdot ,t,x)\right\| _{\mathcal{H}_{1}^{\otimes
n}}^{2}=E\left[ \left( \int_{0}^{t}\delta _{0}(B_{2s})ds\right) ^{n}\right]
\text{. }
\end{equation*}

To handle the case  $H>\frac{1}{2}$, we need the following technical lemma.

\begin{lemma}
\label{lem1} Set%
\begin{equation}
g_{\mathbf{s}}(x_{1},\ldots ,x_{n})=p_{t-s_{\sigma (n)}}(x-x_{\sigma
(n)})\cdots p_{s_{\sigma (2)}-s_{\sigma (1)}}(x_{\sigma (2)}-x_{\sigma
(1)})).  \label{ge}
\end{equation}
Then,%
\begin{equation*}
\left\langle g_{\mathbf{s}},g_{\mathbf{t}}\right\rangle _{L^{2}(\mathbb{R}%
^{nd})}=\psi (\mathbf{s},\mathbf{t}),
\end{equation*}%
where $\psi (\mathbf{s},\mathbf{t})$ is defined in (\ref{a1c}).
\end{lemma}

\begin{proof}
\bigskip By Plancherel's identity%
\begin{equation*}
\left\langle g_{\mathbf{s}},g_{\mathbf{t}}\right\rangle _{L^{2}(\mathbb{R}%
^{nd})}=(2\pi )^{-dn}\left\langle \mathcal{F}g_{\mathbf{s}},\mathcal{F}g_{%
\mathbf{t}}\right\rangle _{L^{2}(\mathbb{R}^{nd})},
\end{equation*}%
where $\mathcal{F}$ denotes the Fourier transform, given by%
\begin{eqnarray*}
\mathcal{F}g_{\mathbf{s}}(\xi _{1},\ldots ,\xi _{n}) &=&(2\pi )^{-\frac{nd}{2%
}}\prod_{j=1}^{n}(s_{\sigma (j+1)}-s_{\sigma (j)})^{-\frac{d}{2}} \\
&&\times \int_{\mathbb{R}^{nd}}\prod_{j=1}^{n}\exp \left( i\left\langle \xi
_{j},x_{j}\right\rangle -\frac{\left| x_{\sigma (j+1)}-x_{\sigma (j)}\right|
^{2}}{2\left( s_{\sigma (j+1)}-s_{\sigma (j)}\right) }\right) d\mathbf{x,}
\end{eqnarray*}%
with the convention $x_{n+1}=x$ and $s_{n+1}=t$. Making the change of
variables $u_{j}=x_{\sigma (j+1)}-x_{\sigma (j)}$ if $1\leq j\leq n-1$, and $%
u_{n}=x-x_{\sigma (n)}$, we obtain%
\begin{eqnarray*}
\mathcal{F}g_{\mathbf{s}}(\xi _{1},\ldots ,\xi _{n}) &=&(2\pi )^{-\frac{nd}{2%
}}\prod_{j=0}^{n}(s_{\sigma (j+1)}-s_{\sigma (j)})^{-\frac{d}{2}} \\
&& \!\!\!\!\!\!\!  \!\!\!\!\!\!\!  \!\!\!\!\!\!\!  \!\!\!\!\!\!\!
 \times \int_{\mathbb{R}^{nd}}\prod_{j=1}^{n}\exp \left( i\left\langle \xi
_{\sigma (j)},x-u_{n}-\cdots -u_{j}\right\rangle -\frac{\left| u_{j}\right|
^{2}}{2\left( s_{\sigma (j+1)}-s_{\sigma (j)}\right) }\right) d\mathbf{u} \\
&=&E\left( \prod_{j=1}^{n}\exp \left( i\left\langle \xi _{\sigma
(j)},x-B_{t}-B_{s_{\sigma (j)}}\right\rangle \right) \right)  \\
&=&E\left( \prod_{j=1}^{n}\exp \left( i\left\langle \xi
_{j},x-B_{t}-B_{s_{j}}\right\rangle \right) \right) .
\end{eqnarray*}%
As a consequence,
\begin{equation*}
\left\langle g_{\mathbf{s}},g_{\mathbf{t}}\right\rangle _{L^{2}(\mathbb{R}%
^{nd})}=(2\pi )^{-nd}\int_{\mathbb{R}^{nd}}E\left( \prod_{j=1}^{n}\exp
\left( i\left\langle \xi _{j},B_{s_{j}}^{1}-B_{t_{j}}^{2}\right\rangle
\right) \right) d\mathbf{\xi ,}
\end{equation*}%
which implies the desired result.
\end{proof}

In the case $H>\frac{1}{2}$,\ and assuming that $u_{0}=1$, the next
proposition shows that the norm $(n!)^{2}\left\| f_{n}(\cdot ,t,x)\right\| _{%
\mathcal{H}_{d}^{\otimes n}}^{2}$ \ coincides with the $n$th moment of the
intersection local time of two independent $d$-dimensional Brownian motions
with weight $\phi (t,s)$.

\begin{proposition}
Suppose that $H>\frac{1}{2}$ and $d<4H$. Then, for all $n\geq 1$%
\begin{equation}
(n!)^{2}\left\| f_{n}(\cdot ,t,x)\right\| _{\mathcal{H}_{d}^{\otimes
n}}^{2}\leq \left\| u_{0}\right\| _{\infty }^{2}E\left[ \left(
\int_{0}^{t}\int_{0}^{t}\phi (s,r)\delta
_{0}(B_{s}^{1}-B_{r}^{2})dsdr\right) ^{n}\right] <\infty ,  \label{es}
\end{equation}%
with equality if $u_{0}$ is constant. Moreover, we have:

\begin{enumerate}
\item If $d=1$, there exists a unique solution to Equation (\ref{e2}). \

\item If $d=2$ , then there exists a unique solution in an interval $[0,T]$
provided $T<T_{0}$, where%
\begin{equation}
T_{0}=\left( \beta _{H}\Gamma (1-\frac{d}{4H})^{2H}\right) ^{-1/(2H-1)}.
\end{equation}
\end{enumerate}
\end{proposition}

\begin{proof}
We have
\begin{equation}
(n!)^{2}\left\| f_{n}(\cdot ,t,x)\right\| _{\mathcal{H}_{d}^{\otimes
n}}^{2}\leq \left\| u_{0}\right\| _{\infty
}^{2}\int_{[0,t]^{n}}\prod_{j=1}^{k}\phi (s_{j},t_{j})\left\langle g_{%
\mathbf{s}},g_{\mathbf{t}}\right\rangle _{L^{2}(\mathbb{R}^{nd})}d\mathbf{s}d%
\mathbf{t,}  \label{f1}
\end{equation}%
where $g_{\mathbf{s}}$ is defined in (\ref{ge}).   Then the results
follow easily from   from Lemma \ref{lem1} and Proposition
\ref{prop5}.
\end{proof}

In the two-dimensional case and assuming $H>\frac{1}{2}$,   the
solution would exists in  any interval $[0,T]$ as a distribution in the
Watanabe space \ $\mathbb{D}^{\alpha ,2}$ for any $\alpha >0$ (see \cite{Wat}%
).

\subsection{Case $H<\frac{1}{2}$ and $d=1$}

We know that in this case, the norm in the space $\mathcal{H}$ is defined in
terms of fractional derivatives. The aim of this section is to show that
$\left\| f_{n}(\cdot ,t,x)\right\| _{\mathcal{H}_{1}^{\otimes n}}^{2}$ is
related to  the $n$th moment of a fractional derivative of the
self-intersection local time of two independent one-dimensional Brownian
motions,  and these moments are finite for all $n\geq 1$, provided
  $\frac{3}{8}<H<\frac{1}{2}$.

Consider the operator $\left( K_{H}^{\ast }\right) ^{\otimes 2}$  on
functions of two variables defined as the action of the operator $%
K_{H}^{\ast }$ on each coordinate. That is, using the notation (\ref{ke3})
we have
\begin{eqnarray*}
&&\left( K_{H}^{\ast }\right) ^{\otimes
2}f(r_{1},r_{2})=K_{H}(T,r_{1})K_{H}(T,r_{2})f(r_{1},r_{2}) \\
&&+K_{H}(T,r_{1})\int_{r_{2}}^{t}\frac{\partial K_{H}}{\partial s}%
(s,r_{2})\left( f(r_{1},s)-f(r_{1},r_{2})\right) ds \\
&&+K_{H}(T,r_{2})\int_{r_{1}}^{t}\frac{\partial K_{H}}{\partial s}%
(v,r_{1})\left( f(v,r_{2})-f(r_{1},r_{2})\right) dv \\
&&+\int_{r_{2}}^{t}\int_{r_{1}}^{t}\frac{\partial K_{H}}{\partial s}(s,r_{2})%
\frac{\partial K_{H}}{\partial v}(v,r_{1})\left[
f(v,s)-f(r_{1},s)-f(v,r_{2})-f(r_{1},r_{2})\right] dsdv.
\end{eqnarray*}%
Suppose that $f(s,t)$ is a continuous function on $[0,T]^{2}$. Define the H%
\"{o}lder norms
\begin{equation*}
\left\| f\right\| _{1,\gamma }=\sup \left\{ \frac{\left|
f(s_{1},t)-f(s_{2},t)\right| }{|s_{1}-s_{2}|^{\gamma }},s_{1},s_{2},t\in
T,s_{1}\neq s_{2}\right\} ,
\end{equation*}%
\begin{equation*}
\left\| f\right\| _{2,\gamma }=\sup \left\{ \frac{\left|
f(s,t_{1})-f(s,t_{2})\right| }{|t_{1}-t_{2}|^{\gamma }},t_{1},t_{2},s\in
T,t_{1}\neq t_{2}\right\}
\end{equation*}%
and%
\begin{eqnarray*}
\left\| f\right\| _{1,2,\gamma } &=&\sup \frac{\left|
f(s_{1},t_{1})-f(s_{1},t_{2})-f(s_{2},t_{1})+f(s_{2},t_{2})\right| }{%
|s_{1}-s_{2}|^{\gamma }|t_{1}-t_{2}|^{\gamma }},  \\
\end{eqnarray*}%
where the supremum is taken in the set $\left\{ t_{1},t_{2},s_{2},s_{2}\in
T,s_{1}\neq s_{2},t_{1}\neq t_{2}\right\} $. Set%
\begin{equation*}
\left\| f\right\| _{0,\gamma }=\left\| f\right\| _{1,\gamma }+\left\|
f\right\| _{2,\gamma }+\left\| f\right\| _{1,2,\gamma }
\end{equation*}%
Then, $\left( K_{H}^{\ast }\right) ^{\otimes 2}f$ is well defined if \ $%
\left\| f\right\| _{0,\gamma }<\infty $ for some $\gamma >\frac{1}{2}-H$. As
a consequence, if $B^{1}$ and $B^{2}$ are two independent one-dimensional
Brownian motions, the following random variable is well defined for all $%
\varepsilon >0$%
\begin{equation}
J_{\varepsilon }=\int_{0}^{T}\left( K_{H}^{\ast }\right) ^{\otimes
2}p_{\varepsilon }(B_{\cdot }^{1}-B_{\cdot }^{2})(r,r)dr.  \label{v7}
\end{equation}%
The next theorem asserts that $J_{\varepsilon }$ converges in $L^{p}$ for
all $p\ge 2$ to a fractional derivative of the intersection local time of \ $B^{1}
$ and $B^{2}$.

\begin{proposition}
\label{prop1} Suppose that $\frac{3}{8}<H<\frac{1}{2}$.Then, for any integer
$k\geq 1$ and, $T>0$ we have $E\left( J_{\varepsilon }^{k}\right) \geq 0$
and
\begin{equation*}
\sup_{\varepsilon >0}E\left( J_{\varepsilon }^{k}\right) <\infty .
\end{equation*}%
Moreover, for all $p\geq 2$, $J_{\varepsilon }$ converges in $L^{p}$ as $%
\varepsilon $ tends to zero to a random variable denoted by%
\begin{equation*}
\int_{0}^{T}\left( K_{H}^{\ast }\right) ^{\otimes 2}\delta
_{0}(B^{1} _\cdot-B^{2}_\cdot )(r,r)dr.
\end{equation*}
\end{proposition}

\begin{proof}
Fix $k\geq 1$. Let us compute the moment of order $k$ of $J_{\varepsilon }$.
We can write%
\begin{equation}
E\left( J_{\varepsilon }^{k}\right) =\int_{[0,T]^{k}}E\left(
\prod_{i=1}^{k}\left( K_{H}^{\ast }\right) ^{\otimes 2}p_{2\varepsilon
}(B^{1}-B^{2})(r_{i},r_{i})\right) d\mathbf{r}.  \label{y2}
\end{equation}%
Using the expression (\ref{ke1}) for the operator $K_{H}^{\ast }$, and the
notation (\ref{a1c}) yields
\begin{equation}
E\left( J_{\varepsilon }^{k}\right) =\int_{[0,T]^{3k}}\psi _{\varepsilon }(%
\mathbf{s},\mathbf{t})\prod_{i=1}^{k}K_{H}^{\ast }(ds_{i},r_{i})K_{H}^{\ast
}(dt_{i},r_{i})d\mathbf{r}.  \label{y3}
\end{equation}%
As a consequence, using (\ref{z3}) we obtain%
\begin{eqnarray*}
E\left( J_{\varepsilon }^{k}\right)  &=&(2\pi )^{-k}\int_{[0,T]^{k}}\int_{%
\mathbb{R}^{k}}\int_{[0,T]^{2k}}e^{-\sum_{j,l=1}^{k}\xi _{j}\xi _{l}\mathrm{%
Cov}\left( B_{s_{j}}^{1}-B_{t_{j}}^{2},B_{s_{l}}^{1}-B_{t_{l}}^{2}\right) }
\\
&&\times \prod_{i=1}^{k}K_{H}^{\ast }(ds_{i},r_{i})K_{H}^{\ast
}(dt_{i},r_{i})e^{-\frac{\varepsilon }{2}\sum_{j=1}^{k}\xi _{j}^{2}}d\mathbf{%
\xi }d\mathbf{r} \\
&\leq &(2\pi )^{-k}\int_{[0,T]^{k}}\int_{\mathbb{R}^{k}}\left(
\int_{[0,T]^{k}}e^{-\frac{1}{2}\mathrm{Var}\left( \sum_{j=1}^{k}\xi
_{j}B_{t_{j}}^{1}\right) }\prod_{i=1}^{k}K_{H}^{\ast }(dt_{i},r_{i})\right)
^{2}d\mathbf{\xi }d\mathbf{r.}
\end{eqnarray*}%
Then, it suffices to show that for each $k$ the following quantity is finite%
\begin{equation}
\int_{T_{k}}\int_{T_{k}^{2}}\prod_{j=1}^{k}\left[ s_{j}-s_{j-1}+t_{j}-t_{j-1}%
\right] ^{-\frac{1}{2}}\prod_{i=1}^{k}K_{H}^{\ast
}(ds_{i},r_{i})\prod_{i=1}^{k}K_{H}^{\ast }(dt_{i},r_{i})d\mathbf{r,}
\label{m1}
\end{equation}%
where $T_{k}=\{0<t_{1}<\cdots <t_{k}<T\}$. Fix a constant $a>0$. We are
going to compute%
\begin{equation*}
\int_{T_{k}}\prod_{j=1}^{k}\left[ t_{j}-t_{j-1}+a\right] ^{-2}%
\prod_{i=1}^{k}K_{H}^{\ast }(dt_{i},r_{i}).
\end{equation*}%
To do this we need some notation. Let $\Delta _{j}$ \ and $I_{j}$ be the
operators defined \ on a function $f(t_{1},\ldots ,t_{k})$ by%
\begin{equation*}
\Delta _{j}f=f-f|_{t_{j}=r_{j}},
\end{equation*}%
and%
\begin{equation*}
I_{j}f=f|_{t_{j}=r_{j}}.
\end{equation*}%
The operator $K_{H}^{\ast }(dt_{i},r_{i})$ is the sum of two components (see
(\ref{ke3})), and it suffices to consider only the second one because the
first one is easy to control. In this way we need to estimate the following
term
\begin{eqnarray*}
&&\int_{T_{k}}\left[ \int_{[0,T]^{k}}\Delta _{1}\cdots \Delta _{k}\left(
\prod_{j=1}^{k}t_{j}^{H-\frac{1}{2}}\left[ t_{j}-t_{j-1}+a\right] ^{-\frac{1%
}{2}}\mathbf{1}_{\{t_{j-1}<t_{j}\}}\right) \right.  \\
&&\left. \times \prod_{j=1}^{k}(t_{j}-r_{j})^{H-\frac{3}{2}}r_{j}^{\frac{1}{2%
}-H}\mathbf{1}_{\{r_{j}<t_{j}\}}\right] ^{2}d\mathbf{r.}
\end{eqnarray*}%
Because $t_{j}^{H-\frac{1}{2}}r_{j}^{\frac{1}{2}-H}\leq 1$, we can disregard
the  factors $r_{j}^{\frac{1}{2}-H}$ and $t_{j}^{H-\frac{1}{2}}$. Using the
rule
\begin{eqnarray*}
\Delta _{j}(FG) &=&F(t_{j})G(t_{j})-F(r_{j})G(r_{j}) \\
&=&\left[ F(t_{j})-F(r_{j})\right] G(t_{j})+F(r_{j})\left[ G(t_{j})-G(r_{j})%
\right]  \\
&=&\Delta _{j}FG+I_{j}F\Delta _{j}G,
\end{eqnarray*}%
we obtain
\begin{eqnarray*}
&&\Delta _{1}\cdots \Delta _{k}\left( \prod_{i=1}^{k}\left[ t_{j}-t_{j-1}+a%
\right] ^{-\frac{1}{2}}\mathbf{1}_{\{t_{j-1}<t_{j}\}}\right)  \\
&=&\sum_{S}\prod_{j=1}^{k}S_{j}\left( \left[ t_{j}-t_{j-1}+a\right] ^{-\frac{%
1}{2}}\mathbf{1}_{\{t_{j-1}<t_{j}\}}\right) ,
\end{eqnarray*}%
where $S_{j}$ is an operator of the form:

\begin{equation*}
II_{j},I\Delta _{j},\Delta _{j-1}I_{j},\Delta _{j-1}\Delta _{j},
\end{equation*}%
and for each $j$, $\Delta _{j}$ must appear only once in the product
$\prod_{j=1}^{k}S_{j}$. Let us
estimate each one of the possible four terms. Fix $\varepsilon >0$ such that
$H-\frac{3}{8}>2\varepsilon $.

\begin{enumerate}
\item Term $II_{j}$:
\begin{equation*}
II_{j}\left( \left[ t_{j}-t_{j-1}+a\right] ^{-\frac{1}{2}}\mathbf{1}%
_{\{t_{j-1}<t_{j}\}}\right) =\left[ r_{j}-t_{j-1}+a\right] ^{-\frac{1}{2}}%
\mathbf{1}_{\{t_{j-1}<r_{j}\}},
\end{equation*}

\item Term $I\Delta _{j}$: \
\begin{eqnarray*}
&&\left| I\Delta _{j}\left( \left[ t_{j}-t_{j-1}+a\right] ^{-\frac{1}{2}}%
\mathbf{1}_{\{t_{j-1}<t_{j}\}}\right) \right| \\
&=&\left| \left[ t_{j}-t_{j-1}+a\right] ^{-\frac{1}{2}}\mathbf{1}%
_{\{t_{j-1}<t_{j}\}}-\left[ r_{j}-t_{j-1}+a\right] ^{-\frac{1}{2}}\mathbf{1}%
_{\{t_{j-1}<r_{j}\}}\right| \\
&\leq &C\left[ t_{j}-r_{j}\right] ^{\frac{1}{2}-H+\varepsilon }\left[
r_{j}-t_{j-1}+a\right] ^{H-1-\varepsilon }\mathbf{1}_{\{t_{j-1}<r_{j}\}} \\
&&+C\left[ t_{j}-t_{j-1}+a\right] ^{-\frac{1}{2}}\mathbf{1}%
_{\{r_{j}<t_{j-1}\}}.
\end{eqnarray*}

\item Term $\mathbf{\ }\Delta _{j-1}I$:%
\begin{eqnarray*}
&&\left| \Delta _{j-1}I\left( \left[ t_{j}-t_{j-1}+a\right] ^{-\frac{1}{2}}%
\mathbf{1}_{\{t_{j-1}<t_{j}\}}\right) \right| \\
&=&\left| \left[ t_{j}-t_{j-1}+a\right] ^{-\frac{1}{2}}\mathbf{1}%
_{\{t_{j-1}<t_{j}\}}-\left[ t_{j}-r_{j-1}+a\right] ^{-\frac{1}{2}}\mathbf{1}%
_{\{r_{j-1}<t_{j}\}}\right| \\
&\leq &C\left[ t_{j-1}-r_{j-1}\right] ^{\frac{1}{2}-H+\varepsilon }\left[
t_{j}-t_{j-1}+a\right] ^{H-1-\varepsilon }\mathbf{1}_{%
\{r_{j-1}<t_{j-1}<t_{j}\}}.
\end{eqnarray*}

\item Term $\Delta _{j-1}\Delta _{j}$:%
\begin{eqnarray*}
&&\left| \Delta _{j-1}\Delta _{j}\left( \left[ t_{j}-t_{j-1}+a\right] ^{-%
\frac{1}{2}}\mathbf{1}_{\{t_{j-1}<t_{j}\}}\right) \right| \\
&=&\left| \left[ t_{j}-t_{j-1}+a\right] ^{-\frac{1}{2}}\mathbf{1}%
_{\{t_{j-1}<t_{j}\}}-\left[ r_{j}-t_{j-1}+a\right] ^{-\frac{1}{2}}\mathbf{1}%
_{\{t_{j-1}<r_{j}\}}\right. \\
&&\left. -\left[ t_{j}-r_{j-1}+a\right] ^{-\frac{1}{2}}\mathbf{1}%
_{\{r_{j-1}<t_{j}\}}+\left[ r_{j}-r_{j-1}+a\right] ^{-\frac{1}{2}}\mathbf{1}%
_{\{r_{j-1}<r_{j}\}}\right| \\
&\leq &C\left[ t_{j}-r_{j}\right] ^{\frac{1}{2}-H+\varepsilon }\left[
t_{j-1}-r_{j-1}\right] ^{\frac{1}{2}-H+\varepsilon }\left[ r_{j}-t_{j-1}+a%
\right] ^{2H-\frac{3}{2}-2\varepsilon }\mathbf{1}_{\{t_{j-1}<r_{j}<t_{j}\}}
\\
&&+C\left[ t_{j-1}-r_{j-1}\right] ^{\frac{1}{2}-H+\varepsilon }\left[
t_{j}-t_{j-1}+a\right] ^{H-1-\varepsilon }\ \mathbf{1}_{%
\{r_{j}<t_{j-1}<t_{j}\}} \\
&&+C\left[ r_{j}-r_{j-1}+a\right] ^{-\frac{1}{2}}\mathbf{1}%
_{\{r_{j-1}<r_{j}<t_{j-1}<t_{j}\}}.
\end{eqnarray*}
\end{enumerate}

If we replace the constant $a$ by $s_{j}-s_{j-1}$ and we treat the the term
\ $s_{j}-s_{j-1}$ in the same way, using the inequality
\begin{equation*}
(a+b)^{-\alpha }\leq  a^{-\frac{\alpha }{2}}b^{-\frac{\alpha }{2}},
\end{equation*}%
we obtain   the same   estimates as if we had
started with
\begin{equation*}
\int_{T_{k}}\left( \int_{T_{k}}\left[ t_{j}-t_{j-1}\right] ^{-\frac{1}{4}%
}\prod_{j=1}^{k}K_{H}^{\ast }(dt_{j},r_{j})\right) ^{2}d\mathbf{r,}
\end{equation*}%
instead of (\ref{m1}). As a consequence, it suffices to control the
following integral%
\begin{equation}
\int_{T_{k}}\left( \int_{T_{k}}\prod_{j=1}^{k}A_{j}^{a,b}(\mathbf{t,r})d%
\mathbf{t}\right) ^{2}d\mathbf{r,}  \label{m5}
\end{equation}%
where $a,b\in \{0,1\}$, and $A_{j}$ has one of the following forms
\begin{eqnarray*}
A_{j}^{0,0} &=&\left[ r_{j}-t_{j-1}\right] ^{-\frac{1}{4}}\mathbf{1}%
_{\{t_{j-1}<r_{j}\}}, \\
A_{j,1}^{0,1} &=&\left[ t_{j}-r_{j}\right] ^{-1+\varepsilon }\left[
r_{j}-t_{j-1}\right] ^{H-\frac{3}{4}-\varepsilon }\mathbf{1}%
_{\{t_{j-1}<r_{j}\}} \\
A_{j,2}^{0,1} &=&\left[ t_{j}-t_{j-1}\right] ^{-\frac{1}{4}}\left[
t_{j}-r_{j}\right] ^{H-\frac{3}{2}}\mathbf{1}_{\{r_{j}<t_{j-1}\}}, \\
A_{j}^{1,0} &=&\left[ t_{j-1}-r_{j-1}\right] ^{-1+\varepsilon }\left[
t_{j}-t_{j-1}\right] ^{H-\frac{3}{4}-\varepsilon }\mathbf{1}%
_{\{r_{j-1}<t_{j-1}<t_{j}\}}, \\
A_{j,1}^{1,1} &=&\left[ t_{j}-r_{j}\right] ^{-1+\varepsilon }\left[
t_{j-1}-r_{j-1}\right] ^{-1+\varepsilon }\left[ r_{j}-t_{j-1}\right] ^{2H-%
\frac{5}{4}-2\varepsilon }\mathbf{1}_{\{t_{j-1}<r_{j}<t_{j}\}}, \\
A_{j,2}^{1,1} &=&\left[ t_{j-1}-r_{j-1}\right] ^{-1+\varepsilon }\left[
t_{j}-t_{j-1}\right] ^{H-\frac{5}{4}-\varepsilon }\ \left[ t_{j}-r_{j}\right]
^{H-\frac{3}{2}}\mathbf{1}_{\{r_{j}<t_{j-1}<t_{j}\}}, \\
A_{j,3}^{1,1} &=&\left[ r_{j}-r_{j-1}\right] ^{-\frac{1}{4}}\left[
t_{j}-r_{j}\right] ^{H-\frac{3}{2}}\left[ t_{j-1}-r_{j-1}\right] ^{H-\frac{3%
}{2}}\mathbf{1}_{\{r_{j-1}<r_{j}<t_{j-1}<t_{j}\}},
\end{eqnarray*}%
and with the convention that any term of the form $A_{j}^{0,1}$ or $%
A_{j}^{1,1}$ must be followed by $A_{j}^{0,0}$ or $A_{j}^{0,1}$ and any term
of the form \ $A_{j}^{0,0}\ $\ or $A_{j}^{1,0}$ must be followed by $%
A_{j}^{1,0}$ or $A_{j}^{1,1}$.   It is not difficult to check that the
integral \ (\ref{m5}) is finite. For instance, for a product of the form  $A_{j-1}^{0,0}A_{j,1}^{1,1}$ we get%
\begin{eqnarray*}
&&\int_{\{r_{j-1}<t_{j-1}<r_{j}<t_{j}\}}\left[ r_{j-1}-t_{j-2}\right] ^{-%
\frac{1}{4}}\left[ t_{j-1}-r_{j-1}\right] ^{-1+\varepsilon }\left[
r_{j}-t_{j-1}\right] ^{2H-\frac{5}{4}-2\varepsilon }  \\
&&\times \left[ t_{j}-r_{j}\right]
^{-1+\varepsilon }\ dt_{j-1} \\
&=&\left[ r_{j-1}-t_{j-2}\right] ^{-\frac{1}{4}}[r_{j}-r_{j-1}]^{2H-\frac{5}{%
4}-\varepsilon }\left[ t_{j}-r_{j}\right] ^{-1+\varepsilon },
\end{eqnarray*}%
and the integral in the variable $r_{j}$ of the square of  this expression  will be
finite because $4H-\frac{5}{2}-2\varepsilon >-1$.

So, we have proved that $\sup_{\varepsilon }E(J_{\varepsilon }^{k})<\infty $
for all $k$. Notice that all these moments are positive. It holds that $%
\lim_{\varepsilon ,\delta \downarrow 0}E(J_{\varepsilon }J_{\delta })$
exists, and this implies the convergence in $L^{2}$, and also in \ $L^{p}$,
for all $p\geq 2$.
\end{proof}

On the other hand, if the initial condition of Equation (\ref{01})
$\ $\ is a
constant $K$, then for all $n\geq 1$ we have%
\begin{equation*}
(n!)^{2}\left\| f_{n}(\cdot ,t,x)\right\| _{\mathcal{H}_{1}^{\otimes
n}}^{2}=K^{2}E\left[ \left( \int_{0}^{T}\left( K_{H}^{\ast }\right)
^{\otimes 2}\delta _{0}(B^{1}_\cdot-B^{2}_\cdot)(r,r)dr\right) ^{n}\right] <\infty ,
\end{equation*}%
provided $H\in \left( \frac{3}{8},\frac{1}{2}\right) $. In fact, by Lemma %
\ref{lem1} we have%
\begin{eqnarray*}
(n!)^{2}\left\| f_{n}(\cdot ,t,x)\right\| _{\mathcal{H}_{1}^{\otimes n}}^{2}
&=&K^{2}\int_{[0,t]^{2n}}\ \left\langle g_{\mathbf{s}},g_{\mathbf{t}%
}\right\rangle _{L^{2}(\mathbb{R}^{n})}   \prod_{i=1}^{n}K_{H}^{\ast
}(dt_{i},r_{i})   \\
&&\times \prod_{i=1}^{n}K_{H}^{\ast }(ds_{i},r_{i})d\mathbf{s}d\mathbf{%
t} \\
&=&K^{2}\int_{[0,t]^{2n}}\ \psi (\mathbf{s},\mathbf{t}%
)\prod_{i=1}^{n}K_{H}^{\ast }(dt_{i},r_{i})\prod_{i=1}^{n}K_{H}^{\ast
}(ds_{i},r_{i})d\mathbf{s}d\mathbf{t,}
\end{eqnarray*}
and it suffices to apply the above proposition.

However,  we do not know the
rate of convergence of the sequence $\left\| f_{n}(\cdot ,t,x)\right\| _{%
\mathcal{H}_{1}^{\otimes n}}^{2}$ as $n$ tends to infinity, and for this reason we are not able to
show the existence of a solution to Equation (\ref{01}) in this case.

\setcounter{equation}{0}
\section{Moments of the solution}

In this section we introduce an approximation of the Gaussian noise $W^{H}$
by means of an approximation of the identity. In the space variable we
choose the heat kernel to define this approximation and in the time variable
we choose a rectangular kernel. In this way, for any $\varepsilon >0$ and $%
\delta >0$ we set
\begin{equation}
\dot{W}_{t,x}^{{\varepsilon ,\delta }}=\int_{0}^{t}\int_{\mathbb{R} ^d}\varphi
_{{\delta }}(t-s)p_{{\varepsilon }}(x-y)dW_{s,y}^{H},  \label{e3a}
\end{equation}%
where%
\begin{equation*}
\varphi _{{\delta }}(t)=\frac{1}{{\delta }}\mathbf{1}_{[0,{\delta }]}(t).
\end{equation*}%
Now we consider the approximation of Equation (\ref{01}) defined by
\begin{equation}
\frac{\partial u_{t,x}^{{\varepsilon ,\delta }}}{\partial t}=\frac{1}{2}%
\Delta u_{t,x}^{{\varepsilon ,\delta }}+u_{t,x}^{{\varepsilon ,\delta }%
}\diamond \dot{W}_{t,x}^{{\varepsilon ,\delta }}.  \label{e4}
\end{equation}%
We recall that the Wick product $u_{t,x}^{{\varepsilon ,\delta }}\diamond
\dot{W}_{t,x}^{{\varepsilon ,\delta }}$ is well defined as a \ square
integrable random variable provided the random variable $u_{t,x}^{{%
\varepsilon ,\delta }}$ belongs to the space $\mathbb{D}^{1,2}$ (see (\ref%
{z9})), and in this case we have%
\begin{equation}
u_{s,y}^{\varepsilon ,\delta }\diamond \dot{W}_{s,y}^{{\varepsilon ,\delta }%
}=\int_{0}^{s}\int_{\mathbb{R}^d}\varphi _{\delta }(s-r)p_{{\varepsilon }%
}(y-z)u_{s,y}^{\varepsilon ,\delta }\delta W_{r,z}^{H}.  \label{e8}
\end{equation}%
The mild or evolution version of Equation (\ref{e4}) will be%
\begin{equation}
u_{t,x}^{\varepsilon ,\delta }=p_{t}u_{0}(y)+\int_{0}^{t}\int_{\mathbb{R}^d%
}p_{t-s}(x-y)u_{s,y}^{\varepsilon ,\delta }\diamond \dot{W}_{s,y}^{{%
\varepsilon ,\delta }}dsdy.  \label{e8a}
\end{equation}%
Substituting (\ref{e8}) into (\ref{e8a}), and formally applying Fubini's theorem
yields%
\begin{equation}
u_{t,x}^{\varepsilon ,\delta }=p_{t}u_{0}(y)+\int_{0}^{t}\int_{\mathbb{R}^d%
}\left( \int_{0}^{t}\int_{\mathbb{R}^d}p_{t-s}(x-y)\varphi _{\delta }(s-r)p_{{%
\varepsilon }}(y-z)u_{s,y}^{\varepsilon ,\delta }dsdy\right) \delta
W_{r,z}^{H}.  \label{e9a}
\end{equation}
This leads to the following definition.

\begin{definition}
An adapted random field $u^{\varepsilon ,\delta }=\{u_{t,x}^{{\varepsilon
,\delta }},t\geq 0,x\in \mathbb{R}^{d}\}$ is a mild solution to Equation
(\ref{e4}) if   for each $(r,z)\in \mathbb{R}_+\times \mathbb{R}^d$
 the integral
\begin{equation*}
Y_{r,z}^{t,x}= \int_{0}^{t}\int_{\mathbb{R}}p_{t-s}(x-y)\varphi _{\delta }(s-r)p_{{%
\varepsilon }}(y-z)u_{s,y}^{\varepsilon ,\delta }dsdy
\end{equation*}%
exists and  $Y^{t,x}$ is   a Skorohod integrable process such that
 (\ref{e9a}) holds for each $(t,x)$.
\end{definition}

The above definition is equivalent to saying that $u_{t,x}^{\varepsilon
,\delta }\in L^{2}(\Omega )$, and for any random variable $F\in \mathbb{D}%
^{1,2}$ , we have%
\begin{eqnarray*}
E(Fu_{t,x}^{\varepsilon ,\delta }) &=&\ E(F)p_{t}u_{0}(y) \\
&&+\left\langle \left( \int_{0}^{t}\int_{\mathbb{R}^d}p_{t-s}(x-y)\varphi
_{\delta }(s-\cdot )p_{{\varepsilon }}(y-\cdot )u_{s,y}^{\varepsilon ,\delta
}dsdy\right) ,DF\right\rangle _{\mathcal{H}_{d}}.\ \
\end{eqnarray*}

Our aim is to construct a solution of Equation (\ref{e4}) using a suitable
version of Feynman-Kac's formula. Suppose that $B=\left\{ B_{t},t\geq
0\right\} $ is a $d$-dimensional  Brownian motion starting at $0$, independent of $W$%
. Set%
\begin{eqnarray*}
\int_{0}^{t}\dot{W}_{t-s,x+B_{s}}^{\varepsilon ,\delta }ds
&=&\int_{0}^{t}\int_{0}^{t}\int_{\mathbb{R}^d}\varphi _{\delta
}(t-s-r)p_{\varepsilon }(B_{s}+x-y)dW_{r,y}^{H}ds \\
&=&\int_{0}^{t}\int_{\mathbb{R}^d}A_{r,y}^{\varepsilon ,\delta }dW_{r,y}^{H},
\end{eqnarray*}%
where
\begin{equation}
A_{r,y}^{\varepsilon ,\delta }=\int_{0}^{t}\varphi _{\delta
}(t-s-r)p_{\varepsilon }(B_{s}+x-y)ds.  \label{n1}
\end{equation}%
Define
\begin{equation}
u_{t,x}^{\varepsilon ,\delta }=E ^{B}\left( u_{0}(x+B_{t})\exp \left(
\int_{0}^{t}\int_{\mathbb{R}^d}A_{r,y}^{\varepsilon ,\delta }dW_{r,y}^{H}-%
\frac{1}{2}\alpha ^{\varepsilon ,\delta }\right) \right) ,  \label{a51}
\end{equation}%
where $\alpha ^{\varepsilon ,\delta }=\left\| A^{\varepsilon ,\delta
}\right\| _{\mathcal{H}_{d}}^{2}$.

\begin{proposition}
The random field $u_{t,x}^{\varepsilon ,\delta }$ given by (\ref{a51}) is a
solution to Equation (\ref{e4}).
\end{proposition}

\begin{proof}
The proof is based on the notion of $S$ transform from white noise analysis
(see \cite{white}). For any element $\varphi \in \mathcal{H}_{1}$ we define%
\begin{equation*}
S_{t,x}(\varphi )=E \left( u_{t,x}^{\varepsilon ,\delta }F_{\varphi
}\right) ,
\end{equation*}%
where
\begin{equation*}
F_{\varphi }=\exp \left( W^{H}(\varphi )-\frac{1}{2}\left\| \varphi \right\|
_{\mathcal{H}_{d}}^{2}\right) .
\end{equation*}%
From (\ref{a51}) we have%
\begin{eqnarray*}
S_{t,x}(\varphi ) &=&E \left( u_{0}(x+B_{t})\exp \left(
W^{H}(A^{\varepsilon ,\delta }+\varphi )-\frac{1}{2}\alpha ^{\varepsilon
,\delta }-\frac{1}{2}\left\| \varphi \right\| _{\mathcal{H}_{d}}^{2}\right)
\right)  \\
&=&E \left( u_{0}(x+B_{t})\exp \left( \left\langle A^{\varepsilon ,\delta
},\varphi \right\rangle _{\mathcal{H}_{d}}\right) \right)  \\
&=&E  \left( u_{0}(x+B_{t})\exp \left( \int_{0}^{t}\left\langle \varphi
_{\delta }(t-s-\cdot )p_{\varepsilon }(B_{s}+x-\cdot ),\varphi \right\rangle
_{\mathcal{H}_{d}}ds\right) \right) .
\end{eqnarray*}%
By the classical Feynman-Kac's formula, $S_{t,x}(\varphi )$ satisfies the
heat equation with potential $V(t,x)=\left\langle \varphi _{\delta }(t-\cdot
)p_{\varepsilon }(x-\cdot ),\varphi \right\rangle _{\mathcal{H}_{d}}$, that
is,%
\begin{equation*}
\frac{\partial S_{t,x}(\varphi )}{\partial t}=\frac{1}{2}\Delta
S_{t,x}(\varphi )+S_{t,x}(\varphi )\left\langle \varphi _{\delta }(t-\cdot
)p_{\varepsilon }(x-\cdot ),\varphi \right\rangle _{\mathcal{H}_{d}}.
\end{equation*}%
As a consequence,
\begin{equation*}
S_{t,x}(\varphi )=p_{t}u_{0}(x)+\int_{0}^{t}\int_{\mathbb{R}^d}%
p_{t-s}(x-y)S_{s,y}(\varphi )\left\langle \varphi _{\delta }(s-\cdot
)p_{\varepsilon }(y-\cdot ),\varphi \right\rangle _{\mathcal{H}_{d}} dsdy.
\end{equation*}%
Notice that $DF_{\varphi }=\varphi F_{\varphi }$. Hence, for any exponential
random variable of this form we have%
\begin{equation*}
E(u_{t,x}^{\varepsilon ,\delta }F_{\varphi
})=p_{t}u_{0}(x)+\int_{0}^{t}\int_{\mathbb{R}}p_{t-s}(x-y)E\ \left(
u_{t,x}^{\varepsilon ,\delta }\left\langle \varphi _{\delta }(s-\cdot
)p_{\varepsilon }(y-\cdot ),DF_{\varphi }\right\rangle _{\mathcal{H}%
_{d}}\right) ,
\end{equation*}%
and we conclude by the duality relationship between the Skorohod integral
and the derivative operator.
\end{proof}

 The next theorem says that the random variables $u_{t,x}^{%
\varepsilon ,\delta }$ have moments of all orders, uniformly bounded in $%
\varepsilon $ and $\delta $, and converge to the solution to Equation (\ref{01}) as $\delta$ and $\epsilon$ tend to zero. Moreover, it provides an expression for the
 moments of the solution to Equation (\ref{01}).

\begin{theorem}
\label{t1} Suppose that $H\geq \frac{1}{2}$ and $d=1$. Then,
   for any
integer $k\geq 1$ we have
\begin{equation}
\sup_{\varepsilon ,\delta }E\left[ \left| u_{t,x}^{\varepsilon ,\delta
}\right| ^{k}\right] <\infty ,  \label{p1}
\end{equation}%
and the limit $\lim_{\varepsilon \downarrow 0}\lim_{\delta \downarrow
0}u_{t,x}^{\varepsilon ,\delta }$ exists in $L^{p}$, for all $p\geq 1$, and
it coincides with the solution $u_{t,x}$ of Equation (\ref{01}).
Furthermore,  if  $U_0^B(t,x)=\prod
 _{j=1}^{k}u_{0}(x+B_{t}^{j})$, where $B^j$ are independent
 $d$-dimensional Brownian motions,  we have for any $k\geq 2$
\begin{equation}
E\left[ u_{t,x}^{k}\right] =E ^{B}\left[ U_0^B(t,x) \exp \left( \sum_{i<j}\ \int_{0}^{t}\
\delta _{0}(B_{s}^{i}-B_{s}^{j})ds\right) \right] .  \label{p2a}
\end{equation}%
if $H=\frac 12$, and
\begin{equation}
E\left[ u_{t,x}^{k}\right] =E ^{B}\left[ U_0^B(t,x) \exp \left( \ \sum_{i<j}\
\int_{0}^{t}\int_{0}^{t}\phi (s,r)\delta
_{0}(B_{s}^{i}-B_{r}^{j})dsdr\right) \right] .  \label{p3}
\end{equation}
if  $H>\frac{1}{2}$.

In the case $d=2$, for any integer $k\ge 2$ there exists $t_0(k)>0$  such that
for all $t<t_0(k)$ (\ref{p1}) holds.  If $t<t_0(M)$ for some $M\ge 3$ then the limit $\lim_{\varepsilon \downarrow 0}\lim_{\delta \downarrow
0}u_{t,x}^{\varepsilon ,\delta }$ exists in $L^{p}$  for all $2\le p <M$, and
it coincides with the solution $u_{t,x}$ of Equation (\ref{01}). Moreover,
(\ref{p3}) holds for all $1\le k\le M-1$.
 \end{theorem}

\begin{proof}
Fix an integer $k\geq 2$. Suppose that $B^{i}=\left\{ B_{t}^{i},t\geq
0\right\} $, $i=1,\ldots ,k$ are independent $d$-dimensional standard
Brownian motions starting at $0$, independent of $W^H$. \ Then, using (\ref%
{a51}) we have
\begin{equation*}
E\left[ \left( u_{t,x}^{\varepsilon ,\delta }\right) ^{k}\right] =E\ \left(
\prod\limits_{j=1}^{k}E^{B}\left[ u_{0}(x+B_{t}^{j})\exp \left(
\int_{0}^{t}\int_{\mathbb{R}}A_{r,y}^{\varepsilon ,\delta
,B^{j}}dW_{r,y}^{H}-\frac{1}{2}\alpha ^{\varepsilon ,\delta ,B^{j}}\right) %
\right] \right) ,
\end{equation*}%
where $A_{r,y}^{\varepsilon ,\delta ,B^{j}}$ and $\alpha ^{\varepsilon
,\delta ,B^{j}}$ are computed using the Brownian motion $B^{j}$. Therefore,%
\begin{eqnarray*}
E\left[ \left( u_{t,x}^{\varepsilon ,\delta }\right) ^{k}\right]  &=&E\ ^{B}%
\left[ \ \exp \left( \frac{1}{2}\left\| \sum_{j=1}^{k}A^{\varepsilon ,\delta
,B^{j}}\right\| _{\mathcal{H}_{d}}^{2}-\frac{1}{2}\sum_{j=1}^{k}\alpha
^{\varepsilon ,\delta ,B^{j}}\right) \prod\limits_{j=1}^{k}u_{0}(x+B_{t}^{j})%
\right]  \\
&=&E ^{B}\left[ \ \exp \left( \sum_{i<j}\ \left\langle A^{\varepsilon
,\delta ,B^{i}},A^{\varepsilon ,\delta ,B^{j}}\right\rangle _{\mathcal{H}%
_{d}}\right) \prod\limits_{j=1}^{k}u_{0}(x+B_{t}^{j})\right] .
\end{eqnarray*}%
That is, the correction term $\frac{1}{2}\alpha ^{\varepsilon ,\delta }$ in (%
\ref{a51}) due to the Wick product produces a cancellation of the diagonal
elements in the \ square norm of $\sum_{j=1}^{k}A^{\varepsilon ,\delta
,B^{j}}$. The next step is to compute the scalar product $\left\langle
A^{\varepsilon ,\delta ,B^{i}},A^{\varepsilon ,\delta ,B^{j}}\right\rangle _{%
\mathcal{H}_{d}}$ for $i\neq j$. We consider two cases.

\textbf{Case 1.} Suppose first that $H=\frac{1}{2}$ and $d=1$. In this case
we have%
\begin{eqnarray*}
\left\langle A^{\varepsilon ,\delta ,B^{i}},A^{\varepsilon ,\delta
,B^{j}}\right\rangle _{\mathcal{H}_{1}} &=&\ \int_{\mathbb{R}%
}\int_{0}^{t}\int_{0}^{t}\int_{0}^{t}\varphi _{\delta
}(t-s_{1}-r)p_{\varepsilon }(B_{s_{1}}^{i}+x-y) \\
&&\times \varphi _{\delta }(t-s_{2}-r)p_{\varepsilon
}(B_{s_{2}}^{j}+x-y)ds_{1}ds_{2}drdy \\
&=&\ \int_{0}^{t}\int_{0}^{t}\int_{0}^{t}\varphi _{\delta
}(t-s_{1}-r)\varphi _{\delta }(t-s_{2}-r) \\
&&\times p_{2\varepsilon }(B_{s_{1}}^{i}-B_{s_{2}}^{j})ds_{1}ds_{2}dr.
\end{eqnarray*}%
We have%
\begin{eqnarray*}
&&\int_{0}^{t}\varphi _{\delta }(t-s_{1}-r)\varphi _{\delta }(t-s_{2}-r)dr \\
&=&\delta ^{-2}\left[ \left( t-s_{1}\right) \wedge (t-s_{2})-(t-s_{1}-\delta
)^{+}\vee (t-s_{2}-\delta )^{+}\right] ^{+} \\
&=&\eta _{\delta }(s_{1},s_{2}).
\end{eqnarray*}%
It it easy to check that $\eta _{\delta }$ is a a symmetric function on $%
[0,t]^{2}$ such that for any continuous function $g$ on $[0,t]^{2}$,%
\begin{equation*}
\lim_{\delta \downarrow 0}\int_{0}^{t}\int_{0}^{t}\eta _{\delta
}(s_{1},s_{2})g(s_{1},s_{2})ds_{1}ds_{2}=\int_{0}^{t}g(s,s)ds.
\end{equation*}%
As a consequence the following limit holds almost surely
\begin{equation*}
\lim_{\delta \downarrow 0}\left\langle A^{\varepsilon ,\delta
,B^{i}},A^{\varepsilon ,\delta ,B^{j}}\right\rangle _{\mathcal{H}%
_{1}}=\int_{0}^{t}p_{2\varepsilon }(B_{s}^{i}-B_{s}^{j})ds,
\end{equation*}%
and by the properties of the  local time of the one-dimensional
Brownian motion we obtain that, almost surely.%
\begin{equation*}
\lim_{\varepsilon \downarrow 0}\lim_{\delta \downarrow 0}\left\langle
A^{\varepsilon ,\delta ,B^{i}},A^{\varepsilon ,\delta ,B^{j}}\right\rangle _{%
\mathcal{H}_{1}}=\int_{0}^{t}\delta _{0}(B_{s}^{i}-B_{s}^{j})ds.
\end{equation*}%
The function $\eta _{\delta }$ satisfies%
\begin{equation*}
\sup_{0\leq r\leq t}\int_{0}^{t}\eta _{\delta }(s,r)ds\leq 1,
\end{equation*}%
and, as a consequence, the estimate (\ref{a2}) implies that for all $\lambda
>0$
\begin{equation*}
\sup_{\varepsilon ,\delta }E\ ^{B}\left[ \lambda \exp \ \left\langle
A^{\varepsilon ,\delta ,B^{i}},A^{\varepsilon ,\delta ,B^{j}}\right\rangle _{%
\mathcal{H}_{1}}\right] <\infty .
\end{equation*}%
Hence (\ref{p1}) holds and $\lim_{\varepsilon \downarrow 0}\lim_{\delta
\downarrow 0}u_{t,x}^{\varepsilon ,\delta }:=v_{t,x}$ exists in $L^{p}$, for
all $p\geq 1$. Moreover,  $E(v_{t,x}^{k})$ equals to the right-hand side of
Equation (\ref{p2a}). Finally, Equation (\ref{e9a}) and the duality
relationship (\ref{dua}) imply that for any random variable $F\in \mathbb{D}%
^{1,2}$ \ with zero mean we have
\begin{equation*}
E\left( Fu_{t,x}^{\varepsilon ,\delta }\right) =\ E\left( \left\langle
DF,\int_{0}^{t}\int_{\mathbb{R}}\left( \int_{0}^{t}\int_{\mathbb{R}%
}p_{t-s}(x-y)\varphi _{\delta }(s-\cdot )p_{{\varepsilon }}(y-\cdot
)u_{s,y}^{\varepsilon ,\delta }dsdy\right) \right\rangle _{\mathcal{H}%
_{1}}\right) \ ,
\end{equation*}%
and letting $\delta $ and $\varepsilon $ tend to zero we get
\begin{equation*}
E\left( Fv_{t,x}^{{}}\right) =  E\left( \left\langle DF,\int_{0}^{t}\int_{%
\mathbb{R}}\left( \int_{0}^{t}\int_{\mathbb{R}}p_{t-s}(x-y)\varphi _{\delta
}(s-\cdot )p_{{\varepsilon }}(y-\cdot )v_{s,y}dsdy\right) \right\rangle _{%
\mathcal{H}_{1}}\right) \ ,
\end{equation*}%
which implies that the process $v$ is the solution of Equation (\ref{01}),
and by the uniqueness $v_{t,x}=u_{t,x}$.

\textbf{Case 2.} Consider now the case $H>\frac{1}{2}$ and $d=2$.
We have
\begin{eqnarray*}
\left\langle A^{\varepsilon ,\delta ,B^{i}},A^{\varepsilon ,\delta
,B^{j}}\right\rangle _{\mathcal{H}_{d}} &=&\ \int_{\mathbb{R}^2%
}\int_{0}^{t}\int_{0}^{t}\int_{0}^{t}\int_{0}^{t}\varphi _{\delta
}(t-s_{1}-r_{1})p_{\varepsilon }(B_{s_{1}}^{i}+x-y) \\
&&\times \varphi _{\delta }(t-s_{2}-r_{2})p_{\varepsilon
}(B_{s_{2}}^{j}+x-y)ds_{1}ds_{2}\ \phi (r_{1},r_{2})dr_{1}dr_{2}dy \\
&=&\ \int_{0}^{t}\int_{0}^{t}\int_{0}^{t}\int_{0}^{t}\varphi _{\delta
}(t-s_{1}-r_{1})\varphi _{\delta }(t-s_{2}-r_{2}) \\
&&\times p_{2\varepsilon }(B_{s_{1}}^{i}-B_{s_{2}}^{j})ds_{1}ds_{2}\phi
(r_{1},r_{2})dr_{1}dr_{2}.
\end{eqnarray*}%
This scalar product can be written in the following form%
\begin{equation*}
\left\langle A^{\varepsilon ,\delta ,B^{i}},A^{\varepsilon ,\delta
,B^{j}}\right\rangle _{\mathcal{H}_{d}}=\ \int_{0}^{t}\int_{0}^{t}\eta
_{\delta }(s_{1}-s_{2})p_{2\varepsilon
}(B_{s_{1}}^{i}-B_{s_{2}}^{j})ds_{1}ds_{2},
\end{equation*}%
where%
\begin{equation}
\ \eta _{\delta }(s_{1},s_{2})=\int_{0}^{t}\int_{0}^{t}\varphi _{\delta
}(t-s_{1}-r_{1})\varphi _{\delta }(t-s_{2}-r_{2})\ \phi
(r_{1},r_{2})dr_{1}dr_{2}.  \label{a7a}
\end{equation}%
We claim that there exists a constant $\gamma$ such that%
\begin{equation}
\eta _{\delta }(s_{1},s_{2})\leq  \gamma |s_{1}-s_{2}|^{2H-2}.  \label{a8}
\end{equation}%
In fact, if    $\left| s_{2}-s_{1}\right| =s$ we have%
\begin{eqnarray*}
\eta _{\delta }(s_{1},s_{2}) &\leq &H(2H-1)\delta ^{-2}\int_{s-\delta
}^{s}\int_{s}^{s+\delta }|u-v|^{2H-2}dudv \\
&=&\frac{1}{2\delta ^{2}}\left[ (s+\delta )^{2H}-(s-\delta )^{2H}-2s^{2H}%
\right] .
\end{eqnarray*}%
Then
\begin{equation*}
H\delta ^{-2}\int_{s}^{s+\delta }\left( y^{2H-1}-\left( y-\delta \right)
^{2H-1}\right) dy\leq H\delta ^{2H-2}\leq H2^{2-2H}s^{2H-2},
\end{equation*}%
if $s\leq 2\delta $. On the other hand, if \ $s\geq 2\delta $, we have%
\begin{eqnarray*}
\frac{1}{2\delta ^{2}}\left[ (s+\delta )^{2H}-(s-\delta )^{2H}-2s^{2H}\right]
&\leq &\frac{H}{\delta }\left[ s^{2H-1}-(s-\delta )^{2H-1}\right]  \\
&\leq &H(2H-1)(s-\delta )^{2H-2} \\
&\leq &H(2H-1)2^{2-2H}s^{2H-2}.
\end{eqnarray*}%
It it easy to check that for any continuous function $g$ on $[0,t]^{2}$,%
\begin{equation*}
\lim_{\delta \downarrow 0}\int_{0}^{t}\int_{0}^{t}\eta _{\delta
}(s_{1},s_{2})g(s_{1},s_{2})ds_{1}ds_{2}=\int_{0}^{t}\int_{0}^{t}\ \phi
(s_{1},s_{2})g(s_{1},s_{2})ds_{1}ds_{2}.
\end{equation*}%
As a consequence the following limit holds almost surely%
\begin{equation*}
\lim_{\varepsilon \downarrow 0}\lim_{\delta \downarrow 0}\left\langle
A^{\varepsilon ,\delta ,B^{i}},A^{\varepsilon ,\delta ,B^{j}}\right\rangle _{%
\mathcal{H}_{d}}=\int_{0}^{t}\int_{0}^{t}\ \phi (s_{1},s_{2})\delta
_{0}(B_{s_{1}}^{i}-B_{s_{2}}^{j})ds_{1}ds_{2}.
\end{equation*}%
From (\ref{a8}) and the estimate (\ref{v1}) we get%
\begin{equation}
\sup_{\varepsilon ,\delta }E  ^{B}\left[ \exp \left( \lambda \ \left\langle
A^{\varepsilon ,\delta ,B^{i}},A^{\varepsilon ,\delta ,B^{j}}\right\rangle _{%
\mathcal{H}_{d}}\right) \right] <\infty \text{,}  \label{a9}
\end{equation}%
if $\lambda <\lambda _{0}(t)$, where $\lambda _{0}(t)$ is defined in (\ref{v0}%
) with $gamma_T$ replaced by $\gamma$.

Hence,  for any integer $k\ge 2$, if  $t<t_{0}(k)$, where $\frac{k(k-1)}{2}=\lambda
_{0}(t_{0}(k))$, then  (\ref{p1}) holds because
\begin{equation*}
E\left[ \left( u_{t,x}^{\varepsilon ,\delta }\right) ^{k}\right] \leq
\left\| u_{0}\right\| ^{k}\left( E  ^{B}\left[ \exp \left( \frac{k(k-1)}{2}\
\left\langle A^{\varepsilon ,\delta ,B^{1}},A^{\varepsilon ,\delta
,B^{2}}\right\rangle _{\mathcal{H}_{d}}\right) \right] \right) ^{\frac{2}{%
k(k-1)}}.
\end{equation*}%
Finally, if $t<t_0(M)$ and  $M\ge 3$,
the limit  $\lim_{\varepsilon \downarrow 0}\lim_{\delta
\downarrow 0}u_{t,x}^{\varepsilon ,\delta }:=v_{t,x}$ exists in $L^{p}$, for
all $2 \le p <M$ and it is equal  to the right-hand side of Equation (\ref%
{p3}). As in the case $H=\frac{1}{2}$ we show that $v_{t,x}=u_{t,x}$.
\end{proof}

\setcounter{equation}{0}
\section{Pathwise heat equation}

In this section we consider the one-dimensional stochastic partial
differential equation
\begin{equation}
\frac{\partial u}{\partial t}=\frac{1}{2}\Delta u+u\dot{W}_{t,x}^{H},
\label{b1}
\end{equation}%
where the product between  the solution $u$ and the noise $\dot{W}_{t,x}^{H}$
is now an ordinary product. \ We first introduce a notion of solution using
the Stratonovich integral and a weak  formulation of the mild solution.
Given a random field  $v=\{v_{t,x},t\geq 0,x\in \mathbb{R\}}$ such that $%
\int_{0}^{T}\int_{\mathbb{R}}\left| v_{t,x}\right| dxdt<\infty $ a.s. for
all $T>0$, the Stratonovich integral
\begin{equation*}
\int_{0}^{T}\int_{\mathbb{R}}v_{t,x}dW_{t,x}^{H}
\end{equation*}%
is defined as the following limit in probability if it exists
\begin{equation*}
\lim_{\varepsilon \downarrow 0}\lim_{\delta \downarrow 0}\int_{0}^{T}\int_{%
\mathbb{R}}v_{t,x}\dot{W}_{t,x}^{\varepsilon ,\delta }dxdt,
\end{equation*}%
where $W_{t,x}^{\varepsilon ,\delta }$ is the approximation of the noise $%
W^{H}$ introduced in \ (\ref{e3a}).

\begin{definition}
\label{def2} A random field $u=\{u_{t,x},t\geq 0,x\in \mathbb{R\}}$ is a
weak solution to Equation (\ref{b1}) if for any $C^{\infty }$ function $%
\varphi $ with compact support on $\ \mathbb{R}$, we have%
\begin{equation*}
\int_{\mathbb{R}}u_{t,x}\varphi (x)dx=\int_{\mathbb{R}}u_{0}(x)\varphi
(x)dx+\int_{0}^{t}\int_{\mathbb{R}}u_{s,x}\varphi ^{\prime \prime
}(x)dxds+\int_{0}^{t}\int_{\mathbb{R}}u_{s,x}\varphi (x)dW_{s,x}^{H}.
\end{equation*}
\end{definition}

Consider the approximating stochastic heat equation
\begin{equation}
\frac{\partial u^{\varepsilon ,\delta }}{\partial t}=\frac{1}{2}\Delta
u^{\varepsilon ,\delta }+u^{\varepsilon ,\delta }\dot{W}_{t,x}^{\varepsilon
,\delta }.  \label{b2}
\end{equation}

\begin{theorem}
\label{t2} Suppose that $H>\frac{3}{4}$. For any $p\geq 2$, the limit
\begin{equation*}
\lim_{\varepsilon \downarrow 0}\lim_{\delta \downarrow
0}u_{t,x}^{\varepsilon ,\delta }=u_{t,x}
\end{equation*}%
exists in $L^{p}$, and defines a weak solution to Equation (\ref{b2}) in the
sense of Definition \ref{def2}. Furthermore, for any positive integer $k$
\begin{equation*}
E  \left( u_{t,x}^{k}\right) =E^{B}\left[ U^B_0(t,x)
\exp \left(
\sum_{i,j=1}^{k}\int_{0}^{t}\int_{0}^{t}\phi (s_{1},s_{2})\delta
(B_{s_{1}}^{i}-B_{s_{2}}^{j})ds_{1}ds_{2}\newline
\right) \right],
\end{equation*}
where $U^B_0(t,x)$  has been defined in Theorem (\ref{t1}).
\end{theorem}

\begin{proof}
By Feynman-Kac's formula we can write%
\begin{equation}
u_{t,x}^{\varepsilon ,\delta }=E^{B}\left\{ u_{0}(x+B_{t})\exp \left(
\int_{0}^{t}\int_{\mathbb{R}}A_{r,y}^{\varepsilon ,\delta
}dW_{r,y}^{H}\right) \right\} ,  \label{n5}
\end{equation}%
where $A_{r,y}^{\varepsilon ,\delta }$ has been defined in (\ref{n1}). We
will first show that for all $k\geq 1$%
\begin{equation}
\sup_{\delta ,\varepsilon }E\left[ \left| u_{t,x}^{\varepsilon ,\delta
}\right| ^{k}\right] <\infty .  \label{b4}
\end{equation}%
Suppose that $B^{i}=\left\{ B_{t}^{i},t\geq 0\right\} $, $i=1,\ldots ,k$ are
independent standard Brownian motions starting at $0$, independent of $W^{H}$.%
  Then, we have, as in the proof of Theorem \ref{t1}
\begin{equation*}
E\left( \ \left( u_{t,x}^{\varepsilon ,\delta }\right) ^{k}\right) =E\ ^{B}
\left[ \ \exp \left( \frac{1}{2}\sum_{i,j=1}^{k}\ \left\langle
A^{\varepsilon ,\delta ,B^{i}},A^{\varepsilon ,\delta ,B^{j}}\right\rangle _{%
\mathcal{H}_{1}}\right) U^B_0(t,x)  \right] .
\end{equation*}%
Notice that
\begin{equation*}
\left\langle A^{\varepsilon ,\delta B^{i}},A^{\varepsilon ,\delta
B^{j}}\right\rangle _{\mathcal{H}_{1}}=\ \ \int_{0}^{t}\int_{0}^{t}\eta
_{\delta }(s_{1},s_{2})p_{2\varepsilon
}(B_{s_{1}}^{i}-B_{s_{2}}^{j})ds_{1}ds_{2},
\end{equation*}%
where $\eta _{\delta }(s_{1},s_{2})$ satisfies (\ref{a8}). As a consequence,
the inequalities (\ref{a2}) and (\ref{a5}) imply that for all $\lambda >0$,
and all $i$,$j$ we have
\begin{equation*}
\sup_{\varepsilon ,\delta }E\left( \exp \lambda \left\langle A^{\varepsilon
,\delta B^{i}},A^{\varepsilon ,\delta B^{j}}\right\rangle _{\mathcal{H}%
_{1}}\right) <\infty .
\end{equation*}%
Thus, (\ref{b4}) holds, and
\begin{eqnarray*}
&&\lim_{\varepsilon \downarrow 0}\lim_{\delta \downarrow 0}E\left[ \left(
u_{t,x}^{\varepsilon ,\delta }\right) ^{k}\right]  \\
&=&E^{B}\exp \left[ U^B_0(t,x) \exp \left(
\frac{1}{2}\sum_{i,j=1}^{k}\ \int_{0}^{t}\int_{0}^{t}\phi
(s_{1},s_{2})\delta _{0}(B_{s_{1}}^{i}-B_{s_{2}}^{j})ds_{1}ds_{2}\right) %
\right] .
\end{eqnarray*}%
In a similar way we can show that the limit $\lim_{\varepsilon ,\varepsilon
^{\prime }\downarrow 0}\lim_{\delta ,\delta ^{\prime }\downarrow 0}E\left(
u_{t,x}^{\varepsilon ,\delta }u_{t,x}^{\varepsilon ^{\prime },\delta
^{\prime }}\right) $ exists. Therefore, the iterated limit $%
\lim_{\varepsilon \downarrow 0}\lim_{\delta \downarrow 0}E\left[
u_{t,x}^{\varepsilon ,\delta }\right] $ exists in $L^{2}$.

Finally we need to show that%
\begin{equation*}
\lim_{\varepsilon \downarrow 0}\lim_{\delta \downarrow 0}\left(
\int_{0}^{t}\int_{\mathbb{R}}u_{s,x}\varphi
(x)dW_{s,x}^{H}-\int_{0}^{t}\int_{\mathbb{R}}u_{s,x}^{\varepsilon ,\delta
}\varphi (x)\dot{W}_{s,x}^{H}dsdx\right) =0,
\end{equation*}%
in probability. We know that $\int_{0}^{t}\int_{\mathbb{R}%
}u_{s,x}^{\varepsilon ,\delta }\varphi (x)\dot{W}_{s,x}^{H}dsdx$ converges
in $L^{2}$ to some random variable $G$. \ Hence, if
\begin{equation}
B_{\varepsilon ,\delta }=\int_{0}^{t}\int_{\mathbb{R}}\left(
u_{s,x}^{\varepsilon ,\delta }-u_{s,x}\right) \varphi (x)\dot{W}%
_{s,x}^{H}dsdx  \label{m9}
\end{equation}%
converges in $L^{2}$ to zero, $u_{s,x}\varphi (x)$ will be Stratonovich
integrable and
\[
\int_{0}^{t}\int_{\mathbb{R}}u_{s,x}\varphi (x)\dot{W}%
_{s,x}^{H}dsdx=G.
\]
 The convergence to zero of (\ref{m9}) is done as follows.
First we remark that $B_{\varepsilon ,\delta }=\delta (\phi ^{\varepsilon
,\delta })$, where%
\begin{equation*}
\phi _{r,z}^{\varepsilon ,\delta }=\int_{0}^{t}\int_{\mathbb{R}}\left(
u_{s,x}^{\varepsilon ,\delta }-u_{s,x}\right) \varphi (x)\varphi _{\delta
}(s-r)p_{\varepsilon }(x-z)dsdx.
\end{equation*}%
Then, from the properties of the divergence operator, it suffices to show
that%
\begin{equation}
\lim_{\varepsilon \downarrow 0}\lim_{\delta \downarrow 0}E\left( \left\|
D\phi ^{\varepsilon ,\delta }\right\| _{\mathcal{H}_{1}\otimes \mathcal{H}%
_{1}}^{2}\right) =0.  \label{n2}
\end{equation}%
It is clear that $\lim_{\varepsilon \downarrow 0}\lim_{\delta \downarrow
0}E\left( \left\| \phi ^{\varepsilon ,\delta }\right\| _{\mathcal{H}%
_{1}}^{2}\right) =0$. On the other hand,
\begin{equation*}
D\left( \phi _{r,z}^{\varepsilon ,\delta }\right) =\int_{0}^{t}\int_{\mathbb{%
R}}\left( D\left( u_{s,x}^{\varepsilon ,\delta }\right) -D\left(
u_{s,x}\right) \right) \varphi (x)\varphi _{\delta }(s-r)p_{\varepsilon
}(x-z)dsdx,
\end{equation*}%
and%
\begin{equation*}
D\left( u_{s,x}^{\varepsilon ,\delta }\right) =E\ ^{B}\left\{
u_{0}(x+B_{t})\exp \left( \int_{0}^{t}\int_{\mathbb{R}}A_{s,y}^{\varepsilon
,\delta }dW_{s,y}^{H}\right) A^{\varepsilon ,\delta }\right\} .
\end{equation*}%
Then, as before we can show that%
\begin{eqnarray*}
&&\lim_{\varepsilon ,\varepsilon ^{\prime }\downarrow 0}\lim_{\delta .\delta
^{\prime }\downarrow 0}E\left( \left\langle D\left( u_{s,x}^{\varepsilon
,\delta }\right) ,D\left( u_{s,x}^{\varepsilon ^{\prime },\delta ^{\prime
}}\right) \right\rangle _{\mathcal{H}_{1}}^{2}\right) \\
&=&E\ ^{B}\left[ u_{0}(x+B_{t}^{1})u_{0}(x+B_{t}^{2})\exp \left(
\sum_{i,j=1}^{2}\int_{0}^{t}\int_{0}^{t}\phi (s_{1},s_{2})\delta
_{0}(B_{s_{1}}^{i}-B_{s_{2}}^{j})ds_{1}ds_{2}\right) \right. \\
&&\left. \times \int_{0}^{t}\int_{0}^{t}\phi (s_{1},s_{2})\delta
_{0}(B_{s_{1}}^{1}-B_{s_{2}}^{2})ds_{1}ds_{2}\right] .
\end{eqnarray*}%
This implies that $u_{s,x}^{\varepsilon ,\delta }$ converges in the space $%
\mathbb{D}^{1,2}$ to $u_{s,x}$ as $\delta \downarrow 0$ and $\varepsilon
\downarrow 0$. Actually, the limit is in the norm of the space $\mathbb{D}%
^{1,2}(\mathcal{H}_{1})$. Then, (\ref{n2}) follows easily.
\end{proof}

Since the solution is square integrable it admits a Wiener-It\^o chaos
expansion. The explicit form of the Wiener chaos coefficients are given
below.

\begin{theorem}
The solution to \hbox{(\ref{b1})} is given by
\begin{equation}
u_{t,x}=\sum_{n=0}^{\infty }I_{n}(f_{n}(\cdot ,t,x))  \label{n3}
\end{equation}%
where
\begin{eqnarray}
&&f_{n}(t_{1},x_{1},\ldots ,t_{n},x_{n},t,x)  \notag \\
&=&E^{B}\left[ u_{0}(x+B_{t})\exp \left( \frac{1}{2}\int_{0}^{t}\int_{0}^{t}%
\phi (s_{1},s_{2})\delta _{0}(B_{s_{1}}-B_{s_{2}})ds_{1}ds_{2}\right)
\right.   \notag \\
&&\quad \left. \times \delta _{0}(B_{t_{1}}+x-x_{1})\cdots \delta
_{0}(B_{t_{n}}+x-x_{n})\right] .  \label{n4}
\end{eqnarray}
\end{theorem}

\begin{proof}
From the Feynman-Kac formula it follows that
\begin{eqnarray*}
u_{t,x}^{{\varepsilon },{\delta }} &=&E^{B}\left( u_{0}(x+B_{t})\exp \left(
\int_{\mathbb{R}^{2}}A_{r,y}^{\varepsilon ,\delta }dW_{r,y}^{H}\right)
\right) \\
&=&E^{B}\left\{ u_{0}(x+B_{t})\exp \left( \frac{1}{2}\Vert A^{\varepsilon
,}\Vert _{\mathcal{H}_{1}}^{2}\right) \exp \left( \int_{0}^{t}\int_{\mathbb{R%
}}A_{r,y}^{\varepsilon ,\delta }dW_{r,y}^{H}-\frac{1}{2}\Vert A^{\varepsilon
,\delta }\Vert _{\mathcal{H}_{1}}^{2}\right) \right\} \\
&=&\sum_{n=0}^{\infty }I_{n}(f_{n}^{{\varepsilon },{\delta }}(t,x)),
\end{eqnarray*}%
where
\begin{equation*}
f_{n}^{{\varepsilon },{\delta }}(t_{1},x_{1},\ldots ,t_{n},x_{n},t,x)=E^{B}%
\left[ u_{0}(x+B_{t})\exp \left( \frac{1}{2}\Vert A^{\varepsilon ,\delta
}\Vert _{\mathcal{H}_{1}}^{2}\right) \ A_{t_{1},x_{1}}^{\varepsilon ,\delta
}\cdots A_{t_{n},x_{n}}^{\varepsilon ,\delta }\right] .
\end{equation*}%
Letting ${\delta }$ and ${\varepsilon }$ go to $0$, we obtain the chaos
expansion of $u_{t,x}$.
\end{proof}

Consider the stochastic partial differential equation (\ref{b1}) and its
approximation (\ref{b2}). The initial condition is $u_{0}(x)$. We shall
study the strict positivity of the solution. In particular we shall show
that $E\ \left[ |u_{t}(x)|^{-p}\right] <\infty $.

\begin{theorem}
Let $H>3/4$. If $E  \left( |u_{0}(B_{t})|\right) >0 $, then for any $%
0<p<\infty $, we have that
\begin{equation}
E\left( \ |u_{t,x}|^{-p}\right) <\infty  \label{c1}
\end{equation}%
and moreover,
\begin{eqnarray}
E\ \left[ |u_{t}(x)|^{-p}\right] &\leq &\left( E  |u_{0}(x+B_{t})|\right)
^{-p-1}  \notag
  E\ ^{B}\Bigg[ |u_{0}(x+B_{t})|   \\
  && \!\!\!\!\!\!\! \!\!\!\!\!\!\!  \times \exp \left( \frac{p^{2}}{2}%
\int_{0}^{t}\int_{0}^{t}\delta (B_{s_{1}}-B_{s_{2}})\phi
(s_{1},s_{2})ds_{1}ds_{2}\right) \Bigg] .  \label{c2}
\end{eqnarray}
\end{theorem}

\begin{proof}
Denote $\kappa _{p}=\left( E\ ^{B}\left( \left| u_{0}(x+B_{t})\right|
\right) \right) ^{-p-1}$. Then, Jensen's inequality applied to the equality $%
u_{t,x}^{\varepsilon ,\delta }=E\ ^{B}\left\{ u_{0}(x+B_{t})\exp \left(
\int_{0}^{t}\int_{\mathbb{R}}A_{r,y}^{\varepsilon ,\delta
}dW_{r,y}^{H}\right) \right\} $ implies that
\begin{equation*}
|u_{t,x}^{\varepsilon ,\delta }|^{-p}\leq \kappa _{p}E\ ^{B}\left\{
|u_{0}(x+B_{t})|\exp \left( -p\int_{0}^{t}\int_{\mathbb{R}%
}A_{r,y}^{\varepsilon ,\delta }dW_{r,y}^{H}\right) \right\} .
\end{equation*}%
Therefore
\begin{eqnarray*}
E\left[ |u_{t,x}^{\varepsilon ,\delta }|^{-p}\right] &\leq &\kappa _{p}E\
^{B}\left\{ |u_{0}(x+B_{t})|E\left[ \exp \left( -p\int_{0}^{t}\int_{\mathbb{R%
}}A_{r,y}^{\varepsilon ,\delta }dW_{r,y}^{H}\right) \right] \right\} \\
&=&\kappa _{p}E\ ^{B}\left\{ |u_{0}(x+B_{t})|E\left[ \exp \left( \frac{p^{2}%
}{2}\left\| A^{\varepsilon ,\delta }\right\| _{\mathcal{H}_{1}}^{2}\right) %
\right] \right\} ,
\end{eqnarray*}
and we can conclude as in the proof of \ Theorem \ref{t2}.
\end{proof}

Using the theory of rough path analysis (see \cite{Ly}) and $p$-variation
estimates, Gubinelli, Lejay and Tindel \cite{GLT} have proved that for $H>%
\frac{3}{4}$, the equation%
\begin{equation*}
\frac{\partial u}{\partial t}=\frac{1}{2}\Delta u+\sigma (u)\dot{W}_{t,x}^{H}
\end{equation*}%
had a unique mild solution up to a random explosion time $T>0$, provided $%
\sigma \in C_{b}^{2}(\mathbb{R})$. In this sense, the restriction $H>\frac{3%
}{4}$, that we found in the case $\sigma (x)=x$ is natural, and in this
particular case, using chaos expansion and Feynman-Kac's formula we have
been able to show the existence of a solution for all times.

\bigskip
{\noindent \bf Acknowledgment}. We thank Carl Mueller for
discussions.

\end{document}